\documentclass[12pt,leqno]{article}
\usepackage{amsmath,amssymb,amsthm,amscd}
\usepackage{dsfont}
\numberwithin{equation}{section} \allowdisplaybreaks
\begin{document}
\newtheorem{theorem}{Theorem}[section]
\newtheorem{defin}{Definition}[section]
\newtheorem{prop}{Proposition}[section]
\newtheorem{corol}{Corollary}[section]
\newtheorem{lemma}{Lemma}[section]
\newtheorem{rem}{Remark}[section]
\newtheorem{example}{Example}[section]
\title{Isotropic subbundles of $TM\oplus T^*M$}
\author{{\small by}\vspace{2mm}\\Izu Vaisman}
\date{}
\maketitle
{\def\thefootnote{*}\footnotetext[1]%
{{\it 2000 Mathematics Subject Classification: } 53C99, 53D17.
}}
\begin{center} \begin{minipage}{12cm}
A{\footnotesize BSTRACT. We define integrable, big-isotropic
structures on a manifold $M$ as subbundles $E\subseteq TM\oplus
T^*M$ that are isotropic with respect to the natural, neutral metric
(pairing) $g$ of $TM\oplus T^*M$ and are closed by Courant brackets
(this also implies that $[E,E^{\perp_g}]\subseteq E^{\perp_g}$). We
give the interpretation of such a structure by objects of $M$, we
discuss the local geometry of the structure and we give a reduction
theorem.}
\end{minipage}
\end{center}
\section{Introduction}
All the manifolds and mappings of this paper are assumed of
the $C^\infty$ class and the following general notation is used:
$M$ is an $m$-dimensional manifold, $\chi^k(M)$ is the space of
$k$-vector fields, $\Omega^k(M)$ is the space of differential
$k$-forms, $\Gamma$ indicates the space of global cross sections
of a vector bundle, $X,Y,..$ are either contravariant vectors or
vector fields, $\alpha,\beta,...$ are either covariant vectors or
$1$-forms, $d$ is the exterior differential and $L$ is the Lie derivative.
The Einstein summation convention will be used whenever possible.
If required by the context, a linear space $V$ shall be identified
with anyone of the spaces $V\oplus0$, $0\oplus V$.

The vector bundle $T^{big}M=TM\oplus T^*M$ is called the {\it big
tangent bundle}. It has the natural, non degenerate metric of zero
signature (neutral metric)
\begin{equation}\label{gFinC}
g((X,\alpha),(Y,\beta))=\frac{1}{2}(\alpha(Y)+\beta(X)),
\end{equation} the non degenerate, skew-symmetric $2$-form
\begin{equation}\label{omegainC}
\omega((X,\alpha),(Y,\beta))= \frac{1}{2}(\alpha(Y)-\beta(X))
\end{equation} and the Courant bracket of cross sections
\cite{C}
\begin{equation}\label{crosetinC} [(X,\alpha),(Y,\beta)] = ([X,Y],
L_X\beta-L_Y\alpha+\frac{1}{2}d(\alpha(Y)-\beta(X))).\end{equation}

A maximal, $g$-isotropic subbundle $D\subseteq T^{big}M$ is called
an almost Dirac structure and, if $\Gamma(D)$ is closed by the
Courant bracket, $D$ is an integrable or a Dirac structure.
Then, the triple $(D,pr_{TM},[\,,\,])$ is a Lie algebroid.

The almost Dirac structure $D$ produces the generalized
distribution $\mathcal{D}=pr_{TM}D$ endowed with a leaf-wise
differentiable $2$-form $\varpi$ induced by $\omega|_D$.
Conversely, $D$ may be recovered from the pair
$(\mathcal{D},\varpi)$ by means of the formula
\begin{equation}\label{reconstrincazclasic}
D=\{(X,\alpha)\,/\,X\in \mathcal{D},\,\alpha|_{\mathcal{D}}=
i(X)\varpi\}.\end{equation} Furthermore, by a technical
computation that uses (\ref{reconstrincazclasic}), it follows that
$D$ is a Dirac structure iff $\mathcal{D}$ is integrable and the
form $\varpi$ is closed on the leaves of $\mathcal{D}$. Thus, a
Dirac structure on $M$ is equivalent with a generalized foliation
by presymplectic leaves where the leaf-wise presymplectic form is
such that the subbundle (\ref{reconstrincazclasic}) is
differentiable \cite{C}.

While the Dirac structures were introduced as a framework for constrained
dynamics, it is rather the geometry of these structures and their
integrability to a conveniently equipped Lie groupoid that were the
object of numerous studies.

The aim of the present paper is to understand the geometry of a
$g$-isotropic subbundle $E\subseteq T^{big}M$, where the maximality
requirement is dropped; we call them big-isotropic structures. Then,
the subbundle $E$ must be discussed in conjunction with its
$g$-orthogonal bundle $E'\supseteq E$ and the corresponding objects
on $M$ will be a pair of generalized distributions
$\mathcal{E}=pr_{TM}E\subseteq\mathcal{E}'=pr_{TM}E'$ and bilinear
mappings $\varpi_x:\mathcal{E}_x\times\mathcal{E}'_x\rightarrow
\mathds{R}$ $(\forall x\in M)$ with a skew symmetric restriction to
$\mathcal{E}_x\times\mathcal{E}_x$. Like in the Dirac case, a
big-isotropic structure $E$ will be integrable if $\Gamma E$ is
closed by Courant brackets. From the properties of the Courant
bracket, one can see that if $E$ is integrable $\Gamma E'$ is a
module over the Lie algebra $\Gamma E$.

This definition and the relationship with the triple
$(\mathcal{E},\mathcal{E}',\varpi)$ are made precise in Section 2,
where we also give several examples. In particular, $\mathcal{E}$
is a generalized  foliation and its leaves, called the
characteristic leaves of $E$, inherit a presymplectic form.

In Section 3 we extend the construction of a local, canonical basis of
a Dirac structure given in \cite{DW} to (integrable) big-isotropic structures.

In Section 4 we discuss the pullback of a big-isotropic structure by a
mapping and use this operation and the canonical bases of Section 3 in
order to study the structure induced on a characteristic leaf and that
induced on a local transversal submanifold of the leaf. We define a
property called (strong) local decomposability and extend the Dufour-Wade
proof of the essential uniqueness of the transversal structure of a
characteristic leaf of a Dirac structure to
strongly, locally decomposable, big-isotropic structures.

Finally, in Section 5 we discuss the push-forward of a big-isotropic
structure and conditions that ensure the projectability of a big-isotropic
structure to the space of leaves of a foliation. The results are used in order
to prove a reduction theorem of an integrable, big-isotropic structure
of a manifold $M$ to a structure of a quotient space $N/\mathcal{F}$ of
a submanifold $N\subseteq M$ by a foliation $\mathcal{F}$.
\section{Definitions, examples, first properties}
We generalize the notion of a Dirac structure by giving the
following definition.
\begin{defin}\label{defbigiso} {\rm A $g$-isotropic subbundle
$E\subseteq T^{big}M$ of rank $k$ $(0\leq k\leq m)$ will be called a
{\it big-isotropic structure} on $M$. A big-isotropic structure $E$
is {\it integrable} if $\Gamma E$ is closed by the Courant bracket
operation.}\end{defin}

Let $E'$ be the $g$-orthogonal subbundle $E^{\perp_g}$ of $E$. The
following proposition gives an important property of an integrable,
big-isotropic structure.
\begin{prop}\label{propimp}  If $E$ is an integrable, big-isotropic
structure then, $\forall(X,\alpha)\in\Gamma(E)$,
$\forall(Y,\beta)\in\Gamma(E')$, one has
\begin{equation}\label{crosetEE'}[(X,\alpha),(Y,\beta)]\in\Gamma(E').
\end{equation}\end{prop}
\begin{proof} Consider also $(Z,\gamma)\in\Gamma E$. The Courant
bracket has the following property (axiom (v) of the definition of a
Courant algebroid \cite{LWX})
$$X(g((Y,\beta),(Z,\gamma)))=g([(X,\alpha),(Y,\beta)],(Z,\gamma)) +
g((Y,\beta),[(X,\alpha),(Z,\gamma)]) $$
$$+\frac{1}{2}(Z(g((X,\alpha),(Y,\beta))) +
Y(g((X,\alpha),(Z,\gamma)))).$$ In our case, since $E\perp_g E,
E\perp_g E'$ and $E$ is closed by Courant brackets, the previous
formula reduces to
$$g([(X,\alpha),(Y,\beta)],(Z,\gamma)),$$ hence, (\ref{crosetEE'})
holds.\end{proof}
\begin{rem}\label{obs1ptgraf} {\rm The closure of $E$ with respect to
the Courant bracket is a more complex notion than the closure of a
distribution $\Delta\subseteq TM$ with respect to the Lie bracket.
For instance, in the latter case, any vector field $X\in\Delta$ is
an infinitesimal automorphism of $\Delta$ while, in the former
case, if $(X,\alpha)\in\Gamma E$ then $X$ is an infinitesimal
automorphism of $E$ (i.e., $(L_XY,L_X\beta)\in\Gamma E$, $\forall
(Y,\beta)\in\Gamma E$) iff $d\alpha(Y,Z)=0$ for all $Y\in
pr_{TM}E,Z\in pr_{TM}E'$. This is an easy consequence of the
expression of the Courant bracket and of the isotropy of
$E$.}\end{rem}

\begin{example}\label{exDirac} {\rm For $k=m$, the integrable,
big-isotropic structures are the Dirac structures.}\end{example}
\begin{example}\label{locprod} {\rm Let $M$ be a locally product
manifold with the structural foliations
$\mathcal{F}_1,\mathcal{F}_2$, i.e., each point $x\in M$ has a
neighborhood $U\approx V_1\times V_2$ where $V_a$ is a
neighborhood of $x$ in the leaf of $\mathcal{F}_a$ through $x$
$(a=1,2)$. Equivalently, $M$ has an atlas of local coordinates of
the form $(x^h,y^u)$ such that $ \mathcal{F}_1$ has the local
equations $dy^u=0$ and $ \mathcal{F}_2$ has the local equations
$dx^h=0$. Then $T^{big}M$ is the direct sum of the $g$-orthogonal
subbundles $T\mathcal{F}_a\oplus T^*\mathcal{F}_a$ and, if $E$ is
a maximal isotropic subbundle of $T\mathcal{F}_1\oplus
T^*\mathcal{F}_1$ , $E$ is a big-isotropic structure on $M$ with
the orthogonal subbundle $E'=E\oplus(T\mathcal{F}_2\oplus
T^*\mathcal{F}_2)$. Furthermore, assume $\Gamma E$ has local bases
$(Z_l,\zeta_l)$ such that
\begin{equation}\label{bazaexprod} Z_l=Z_l^h(x)\frac{\partial}{\partial x^h},\;
\zeta_l=\zeta_{lh}(x)dx^h\end{equation} and
that $E$ is Dirac along the leaves of $ \mathcal{F}_1$. Then,
using Definition \ref{defbigiso} and the Courant algebroid
properties of the Courant bracket \cite{{C},{LWX}}, it is easy to
check that the big-isotropic structure $E$ is integrable. For
instance, if $P\in\chi^2(M),\theta\in\Omega^2(M)$ have local
expressions that depend only on $(x^h)$ (which is an invariant
property) and if $[P,P]=0,d\theta=0$, then
$graph(\sharp_P|_{T^*\mathcal{F}_1})$,
$graph(\flat_\theta|_{T\mathcal{F}_1})$ are integrable,
big-isotropic structures of $M$ of the kind described
above.}\end{example}
\begin{example}\label{exemplufol} {\rm For any pair of vector subbundles
$F\subseteq F'\subseteq TM$ , $E=F\oplus ann\,F'$ is a
big-isotropic structure on $M$ with the $g$-orthogonal bundle
$E'=F'\oplus ann\,F$. Furthermore, $E$ is integrable iff $F$ is
tangent to a foliation and $\Gamma F'$ is invariant by Lie
brackets with cross sections of $F$; this means that $F'$ is a
projectable distribution with respect to the foliation $F$, i.e.,
$F'$ is projection-related with distributions of the local spaces
of leaves \cite{VTw}.}\end{example}
\begin{example}\label{exgraphtheta} {\rm Let $S$ be a rank $k$ subbundle
of $TM$ and $\theta\in\Omega^2(M)$ a differential $2$-form. Then,
\begin{equation}\label{Etheta}
E_\theta=graph(\flat_\theta|_{S})=\{(X,\flat_{\theta}X=i(X)\theta)\,/\,X\in S\}
\subseteq T^{big}M\end{equation} is a big-isotropic structure on $M$ with
the $g$-orthogonal bundle
\begin{equation}\label{perpEtheta} E'_\theta=\{(Y,\flat_{\theta}
Y+\gamma)\,/\,Y\in TM,\,\gamma\in ann\,S\}.\end{equation} For the
integrability conditions, we compute the Courant bracket of two
cross sections of $E_\theta$. With the notation of (\ref{Etheta})
and with $X,Y\in S$, we get
\begin{equation}\label{crosetintheta}
[(X,\flat_{\theta} X),(Y,\flat_{\theta} Y)]
=([X,Y],(L_Xi(Y)-L_Yi(X)) \theta\end{equation}
$$+d(\theta(X,Y))) =([X,Y],i([X,Y])\theta+i(X\wedge
Y)d\theta).$$ The final result of (\ref{crosetintheta}) is in
$E_\theta$ iff $S$ is involutive and
\begin{equation}\label{conddiftheta} d\theta(X,Y,Z)=0,\hspace{5mm}
\forall X,Y\in \Gamma S,\forall Z\in\chi^1(M).\end{equation}
Therefore, $E_\theta$ is integrable iff $S$ is a foliation and
$\theta$ satisfies (\ref{conddiftheta}) (in particular, if $\theta$
is closed).}\end{example}
\begin{example}\label{exgraphP} {\rm Let $S^*$ be a subbubdle of rank $k$
of $T^*M$ and $P\in\chi^2(M)$ a differentiable bivector field on $M$. Then
\begin{equation}\label{eqEP} E_P=graph(\sharp_P|_{S^*}) =
\{(\sharp_P\sigma=i(\sigma)P,\sigma)\,/\,\sigma\in S^*\}\end{equation} is
a big-isotropic structure on $M$ with the $g$-orthogonal bundle
\begin{equation}\label{eqE'P} E'_P=\{(\sharp_P\beta+Y,\beta)\,/\,\beta\in T^*M,
Y\in ann\,S^*\}.\end{equation}
For the integrability conditions we recall that $P$ defines the bracket of
$1$-forms:
\begin{equation}\label{croset1forme} \{\alpha,\beta\}_P=
L_{\sharp_P\alpha}\beta-L_{\sharp_P\beta}\alpha
-d(P(\alpha,\beta)),\end{equation} which is related to the Schouten-Nijehuis
bracket $[P,P]$ by the Gelfand-Dorfman formula \cite{Dorf} \begin{equation}\label{GD}
P(\{\alpha,\beta\}_{P},\gamma) =
\gamma([\sharp_P\alpha,\sharp_P\beta])
+\frac{1}{2} [P, P](\alpha,\beta,\gamma).
\end{equation} With (\ref{croset1forme}) and (\ref{GD}) we get the Courant bracket
\begin{equation}\label{CEPE'P} [(\sharp_P\sigma,\sigma),
(\sharp_P\tau,\tau)] = (\sharp_P\{\sigma,\tau\}_P -
\frac{1}{2}i(\sigma\wedge\beta)[P,P],\{\sigma,\beta\}_P),\end{equation}
where $\sigma,\tau\in S^*$. The result of (\ref{CEPE'P} is in
$\Gamma E_P$ iff $S^*$ is closed by the bracket (\ref{croset1forme})
and
\begin{equation}\label{cond2-a} [P,P](\sigma_1,\sigma_2,\beta)=0,
\hspace{3mm}\forall\sigma_1,\sigma_2\in S^*,\forall\beta\in
T^*M.\end{equation} Therefore, $E_P$ is integrable iff $S^*$ is
closed by the $P$-brackets (\ref{croset1forme}) and $P$ satisfies
condition (\ref{cond2-a}) (in particular, $P$ is a Poisson bivector
field).}\end{example}
\begin{example}\label{lift} {\rm The construction indicated in
\cite{Vtg} for the lift of a Dirac structure of a manifold $M$ to
its tangent manifold $TM$ may also be used for a big-isotropic
structure $E,E'$. More exactly, if we look at the locally free
sheaves $\underline{E},\underline{E}'$ of germs of cross sections of
$E,E'$, the formulas
\begin{equation}\label{eqlift}\begin{array}{l}
\underline{tg(E)}=span\{(X^C,\alpha^C),(X^V,\alpha^V)\,/\,(X,\alpha)
\in \underline{E}\},\vspace{2mm}\\
\underline{tg(E)'}=span\{(X^C,\alpha^C),(X^V,\alpha^V)\,/\,(X,\alpha)
\in \underline{E}'\},\end{array}\end{equation} where $C,V$ denote
the complete and vertical lift, respectively, define locally free
sheaves of germs of cross sections of orthogonal subbundles
$tg(E),tg(E)'\subseteq T^{big}(TM)$. The formulas for
scalar products and Courant brackets of lifts established in
\cite{Vtg} show that $tg(E)$ is a big-isotropic
structure on $TM$ with the orthogonal bundle $tg(E)'$ and that, if $E$
is integrable, $tg(E)$ is integrable too.}\end{example}

Now, we shall look for objects of $TM$ that are equivalent with a
big-isotropic structure $E$.

For the algebraic aspects, we refer to a fixed point $x\in M$ and we
associate with $E$ the vector spaces
\begin{equation}\label{distribluiE} \mathcal{E}_x=pr_{T_xM}E_x
\subseteq\mathcal{E}'_x=pr_{T_xM}E'_x.\end{equation}
Then, we define a bilinear mapping $\varpi_x:
\mathcal{E}_x\times\mathcal{E}'_x \rightarrow\mathds{R}$ by means of
the formula
\begin{equation}\label{varpi} \varpi_x(X,Y)=
\omega((X,\alpha),(Y,\beta))=\alpha(Y)=-\beta(X),\end{equation}
where $(X,\alpha)\in E_x,(Y,\beta)\in E'_x$. The last equalities
hold and the result is independent on the choice of $\alpha,\beta$
because $(X,\alpha)\perp_g(Y,\beta)$. Of course,
$\varpi|_{\mathcal{E}_x\times\mathcal{E}_x}$ is skew-symmetric.
Notice that $\varpi_x$ may be identified with a mapping
$\flat_{\varpi_x}:\mathcal{E}_x\rightarrow \mathcal{E}^{'*}_xM$,
which sends $X$ to $i(X)\varpi$ (with an obvious notation) and it
is easy to see that $ker\,\flat_{\varpi_x}=T_xM\cap E_x$. This
implies that, if $\varpi$ is non degenerate (i.e.,
$ker\,\flat_{\varpi_x}=0$) $E_x$ is the graph of a mapping
$pr_{T^*_xM}E_x\rightarrow T_xM$. Similarly, if $E_x$ has the
property $T^*_xM\cap E_x=0$ then $E_x$ is the graph of a mapping
$\mathcal{E}_x\rightarrow T^*_xM$. If we are in one (and the same)
of the two cases above $\forall x\in M$, we will say that $E$ is
of {\it the graph type}.
\begin{example}\label{exptErond} {\rm In Example \ref{locprod} we have
$\mathcal{E}=pr_{TM}E,\mathcal{E}'=\mathcal{E}\oplus
T\mathcal{F}_2$ and $\varpi$ is the extension of the $2$-form of
the almost Dirac structure $E$ along $\mathcal{F}_1$ by the value
$0$ for second arguments in $T\mathcal{F}_2$. In Example
\ref{exemplufol}, $\mathcal{E}=F,\mathcal{E}'=F'$ and
$\varpi(X,Y)=0$. In Example
\ref{exgraphtheta}, $\mathcal{E}=S,\mathcal{E}'=
TM,\varpi=\theta|_{S\times TM}$ and in Example \ref{exgraphP},
$$\mathcal{E}=im(\sharp_P|_{S^*}),  \mathcal{E}'=\mathcal{E}+
ann\,S^*,\varpi(\sharp_P\sigma,\sharp_P\beta+Y)=-P(\sigma,\beta)$$
where $\sigma\in S^*,Y\in ann\,S^*$. In the particular case of
Example \ref{exgraphP} where $P$ is the Lie-Poisson bivector field
of the Lie coalgebra $\mathcal{G}^*$ of the connected Lie group
$G$ (e.g., see \cite{V-carte}) and $S^*=\mathcal{G}'$ is the Lie
subalgebra of the connected subgroup $G'\subseteq G$,
$\mathcal{E}$ are the tangent spaces of the orbits of the
coadjoint action of $G'$ on $\mathcal{G}^*$.}\end{example}
\begin{prop}\label{elemechiv} For any pair of planes
$\mathcal{E}_x\subseteq\mathcal{E}'_x\subseteq T_xM$ and any
bilinear mapping $\varpi_x:\mathcal{E}_x\times\mathcal{E}'_x
\rightarrow\mathds{R}$ with a skew-symmetric restriction to
$\mathcal{E}_x\times\mathcal{E}_x$ there exists a big-isotropic
plane $E_x\subseteq T^{big}_xM$ such that {\rm(\ref{distribluiE}),
(\ref{varpi})} are the
given planes and mapping.\end{prop}
\begin{proof} For the given planes and mapping, put
\begin{equation}\label{eqE} \begin{array}{l}
E_x=\{(X,\alpha)\,/\, X\in\mathcal{E}_x \,\&\,\forall
Y\in\mathcal{E}'_x,\,\alpha(Y)=\varpi_x(X,Y)\}, \vspace{2mm}\\
E'_x=\{(Y,\beta)\,/\, Y\in\mathcal{E}'_x \,\&\,\forall
X\in\mathcal{E}_x,\,\beta(X)=-\varpi_x(X,Y)\}.
\end{array}\end{equation} Obviously, covectors $\alpha,\beta$ as required by (\ref{eqE}) exist, hence, the projection on the first term of a pair defines
epimorphisms $E_x\rightarrow\mathcal{E}_x$,
$E'_x\rightarrow\mathcal{E}'_x$ with the kernels $E_x\cap
T_x^*M=ann\,\mathcal{E}'_x$, $E'_x\cap T_x^*M =ann\,\mathcal{E}_x$,
respectively. Accordingly, one has the exact sequences
\begin{equation}\label{exactseqE}
0\rightarrow ann\,\mathcal{E}'_x\rightarrow E_x\rightarrow
\mathcal{E}_x\rightarrow 0,\hspace{2mm}
0\rightarrow ann\,\mathcal{E}_x\rightarrow E'_x\rightarrow
\mathcal{E}'_x\rightarrow 0,\end{equation}
and we get
\begin{equation}\label{relptdimE}
dim\,E_x=dim\,\mathcal{E}_x+dim\,ann\,\mathcal{E}'_x,\;
dim\,E'_x=dim\,\mathcal{E}'_x+dim\,ann\,\mathcal{E}_x,\end{equation}
which implies
$dim\,E_x+dim\,E'_x=2m.$ Thus, $E'_x,E^{\perp_g}_x$ have the same dimension and, since by (\ref{eqE}) $E'_x\perp_g
E_x$, we get
$E'_x=E^{\perp_g}_x$. Furthermore, the skew-symmetry of
$\varpi$ on $\mathcal{E}_x\times\mathcal{E}_x$ implies
$E_x\subseteq E'_x$, hence, $E_x$ is isotropic. The fact that the given planes and mapping are associated with $E_x,E'_x$ of (\ref{eqE}) via (\ref{distribluiE}), (\ref{varpi}) is obvious.
\end{proof}

Now, starting with the big-isotropic structure $E$, let $x$ vary
on $M$. Since $E,E'$ are differentiable vector bundles, the
generalized, distributions $\mathcal{E},\mathcal{E}'$ defined by
the spaces (\ref{distribluiE}) are differentiable. If
$\mathcal{E}$ is a regular distribution (i.e.,
$dim\,\mathcal{E}_x=const.$, therefore, by (\ref{relptdimE}),
$dim\,\mathcal{E}'_x=const.$ as well), the structure $E$ will be
called {\it regular}. If $E$ is integrable then
$(E,pr_{TM},[\,,\,])$ is a Lie algebroid and, accordingly,
$\mathcal{E}$ is a generalized foliation. Proposition
\ref{propimp} implies that if $X\in\Gamma(\mathcal{E}), Y\in\Gamma
(\mathcal{E}')$ then $[X,Y]\in\Gamma (\mathcal{E}')$.

It is worth formalizing the status of $E'$ as follows since this
may be useful in a discussion of the integrability of $E$ to a Lie
groupoid. Let $B\rightarrow M$ be a vector bundle endowed with an
anchor (morphism) $\rho:B\rightarrow TM$. Assume that there exists
a vector subbundle $A\subseteq B$ endowed with a Lie algebroid
structure $(A,\rho|_A,[\,,\,]_A)$ and there exists an
$\mathds{R}$-bilinear operation $[\,,\,]:\Gamma A\times\Gamma
B\rightarrow\Gamma B$ that reduces to $[\,,\,]_A$ for arguments in
$\Gamma A$. Then, $(B,\rho,[\,,\,])$ will be called a {\it modular
enlargement} of the Lie algebroid $A$ if, $\forall f,h\in
C^\infty(M),a\in\Gamma A, b\in\Gamma B$, the following conditions
are satisfied:\\

$1)$\hspace{1cm}$\rho[a,b]=[\rho a,\rho b],$\vspace{2mm}\\

$2)$\hspace{1cm}$[fa,hb]=fh[a,b]+f((\rho a)h)b-h((\rho b)f)a$,\vspace{2mm}\\

$3)$\hspace{1cm}$[a_1,[a_2,b]]=[[a_1,a_2],b] + [a_2,[a_1,,b]]$\vspace{2mm}\\
(the right hand side of $1)$ is a Lie bracket of vector fields).

With this terminology, if $E$ is an integrable, big-isotropic
structure then $E'$ with $\rho=pr_{TM}$ and the Courant bracket is
a modular enlargement of the Lie algebroid $E$.

We also give the following definition:
\begin{defin}\label{deffolcaract} {\rm The
generalized foliation $ \mathcal{E}$ is the {\it characteristic
foliation} of the integrable, big-isotropic structure $E$ and its
leaves are the {\it characteristic leaves}. The generalized
distribution $ \mathcal{E}'$ is the {\it characteristic module} of
$E$.}\end{defin}
\begin{example}\label{constrdinalgL} {\rm We can extend the construction
of Dirac structures from Lie algebroids \cite{BCWZ} as follows.
Let $A\rightarrow M$ be a Lie algebroid of anchor $\rho_A:A\rightarrow TM$
and bracket $[\,,\,]_A$ and $(B,\rho,[\,,\,])$ a modular enlargement of $A$.
Assume that one also has a {\it co-anchor} $\sigma:B\rightarrow T^*M$ such
that the following properties are satisfied for all $a\in \Gamma A,b\in\Gamma B$:\\

$i)$\hspace{1cm}$<\sigma a,\rho b>=- <\sigma b,\rho a>$,\\

$ii)$\hspace{1cm}$\sigma[a,b]=L_{\rho a}(\sigma b) - L_{\rho b}(\sigma a) +
d<\sigma a,\rho b>$.\vspace{2mm}\\
\noindent Furthermore, assume that, $\forall x\in M$, the morphism $(\rho,\sigma):
B\rightarrow T^{big}M$ satisfies the condition
$$rank(\rho,\sigma)|_{A_x}+rank(\rho,\sigma)|_{B_x}=2m.$$
Then, it is easy to check that $$E_\sigma=\{(\rho a,\sigma a)\,/\,\forall a\in A\}$$
is an integrable, big-isotropic, structure on $M$ with the orthogonal bundle
$$E'_\sigma=\{(\rho b,\sigma b)\,/\,\forall b\in B\}.$$
The characteristic foliation of $E_\sigma$ is $\rho(A)$ and the
characteristic module is
$\rho(B)$. In fact, any integrable, big-isotropic, structure $E$
is of this kind, where
$A=E$, $B=E'$, the brackets are Courant brackets, the anchor is
the projection on $TM$ and the co-anchor is the projection on $T^*M$.}\end{example}

Concerning the mapping $\varpi$, we notice that it has the
following differentiability property: for any characteristic leaf
$S$ of $E$ and for any differentiable vector fields
$X,Y\in\chi^1(M)$, such that $X|_S\in\Gamma\mathcal{E}|_S,
Y|_S\in\Gamma\mathcal{E}'|_S$, $\varpi(X,Y)$ is a differentiable
function on $S$. Indeed, since $\mathcal{E}|_S$, hence
$\mathcal{E}'|_S$ too,  has a constant dimension,
(\ref{exactseqE}) produces exact sequences of differentiable,
vector bundles over $S$. Using differentiable splittings of these
sequences, we see that there are $1$-forms $\alpha,\beta\in
T^*_SM$, which are differentiable along $S$, such that
$(X,\alpha)\in\Gamma E|_S, (Y,\beta)\in\Gamma E'|_S$. Accordingly,
(\ref{varpi}) shows that $\varpi(X,Y)$ is a differentiable
function on $S$. This property will be called {\it leaf-wise
differentiability}. If the structure $E$ is regular, the exact
sequences (\ref{exactseqE}) have differentiable splittings over
$M$, and the functions  $\varpi(X,Y)$ are differentiable on the
whole manifold $M$.

A multilinear mapping
\begin{equation}\label{truncform} \lambda_x:\wedge^{s-1}\mathcal{E}_x
\otimes\mathcal{E}'_x\rightarrow\mathds{R},\end{equation} which is
defined $\forall x\in M$, has a totally skew-symmetric restriction
to $\wedge^{s}\mathcal{E}_x$ and is leaf-wise differentiable with
respect to the characteristic foliation of $E$, will be called a
{\it truncated $s$-form} on $(M,E)$. We will denote by $\Omega^s_{tr}(M,E)$
the space of truncated $s$-forms. Because of the integrability
conditions of $E$, the usual formula for the evaluation of the
exterior differential of a differential form makes sense for
truncated $s$-forms on $(M,E)$ and for arguments in $\chi^1(M)$
that belong to $\mathcal{E}$. Moreover, if $E$ is a regular structure,
the exterior differential also makes sense if the last argument belongs
to $\mathcal{E}'$, while the other arguments belong to $\mathcal{E}$.
Whenever it makes sense, we will denote the differential of a truncated
form by $d_{tr}$, .  In the regular case, $d_{tr}$ is a coboundary morphism
\begin{equation}\label{coboundary} d_{tr}:\Omega^s_{tr}(M,E) \rightarrow
\Omega^{s+1}_{tr}(M,E),\hspace{5mm}d_{tr}^2=0\end{equation} and the
cohomology spaces $H^s_{tr}(M,E)$ of the cochain complex
$(\Omega^s_{tr}(M,E),d_{tr})$ will be the {\it truncated, de Rham
cohomology spaces} of $(M,E)$.
\begin{prop}\label{caractnearly}For any integrable, big-isotropic
structure $E$ one has
\begin{equation}\label{anulpi} d_{tr}\varpi(X_1,X_2,X_3)=0,
\hspace{5mm}\forall X_a\in\mathcal{E},\,
a=1,2,3.\end{equation}  If either $E$ is an almost Dirac structure
or $E$ is a regular big-isotropic structure, $E$  is integrable
iff $\mathcal{E}$ is involutive, the corresponding distribution
$\mathcal{E}'$ is invariant by Lie brackets with vector fields of
$ \mathcal{E}$ and {\rm(\ref{anulpi})} with  $X_3$ replaced by
$Y\in\mathcal{E}'$ holds.\end{prop}
\begin{proof} For almost Dirac structures the result is known \cite{C}.
We prove the result for regular, big-isotropic structures and we
will get the first assertion on the way. It was already shown that
$\mathcal{E}$ and $\mathcal{E}'$ satisfy the required conditions
for any integrable, big-isotropic structure. Now, let
$X_1,X_2\in\mathcal{E},Y\in\mathcal{E}'$ be differentiable vector
fields on $M$ and $\alpha_1,\alpha_2,\beta$ differentiable
$1$-forms such that $(X_1,\alpha_1),(X_2,\alpha_2)\in E,
(Y,\beta)\in E'$. (The existence of $\alpha,\beta$ is ensured by
the regularity hypothesis.) Then, keeping in mind the
$g$-orthogonality relations among these pairs and using
(\ref{varpi}), we get
$$(d_{tr}\varpi)(X_1,X_2,Y)=
(d_C\omega)((X_1,\alpha_1),(X_2,\alpha_2),(Y,\beta))$$
$$=X_1(\alpha_2(Y))-X_2(\alpha_1(Y)) + Y(\alpha_1(X_2))
+\alpha_1([X_2,Y]) -\alpha_2([X_1,Y])$$
$$-<L_{X_1}\alpha_2-L_{X_2}\alpha_1-d(\alpha_1(X_2)),Y>=0,$$
where $d_C$ is the  operator defined by the usual expression of
the exterior differential of a Lie algebroid with Courant brackets
instead of Lie brackets.

The previous calculation makes sense and remains true in the
non regular case if we also assume $Y\in\mathcal{E}$. This proves
the first assertion of the proposition.

Back to the regular case, if we start with a big-isotropic structure $E$ which
satisfies the required hypotheses then, for arguments as above, we
get
$$d_{tr}\varpi(X_1,X_2,Y)=2g([(X_1,\alpha_1),
(X_2,\alpha_2)],(Y,\beta))$$ $$=2g([(X_2,\alpha_2),(Y,\beta)],
(X_1,\alpha_1)).$$ Thus, if $\varpi$ is $d_{tr}$-closed, $E$ is closed by
Courant brackets and $E'$ is closed by Courant brackets with cross
sections of $E$.

For the non regular case, we do not get the converse result; from
condition (\ref{anulpi})  it only follows that $[\Gamma E,\Gamma
E]
\in\Gamma E'$.\end{proof}
\begin{corol}\label{Diracasoc} For any big-isotropic structure $E$
and for each point $x\in M$, there exists a canonical extension of
$E_x$ to an almost Dirac space $D_x(E)\subseteq T^{big}_xM$. If
differentiable with respect to $x$, these spaces define a
canonical almost Dirac extension $D(E)$ of $E$ and if $E$ is
integrable so is $D(E)$.\end{corol}
\begin{proof} The restriction of $\varpi$ to $\mathcal{E}\times\mathcal{E}$
is a leaf-wise differentiable $2$-form and we may use (\ref{reconstrincazclasic})
and  define
\begin{equation}\label{eqDE} D_x(E)=\{(X,\alpha)\,/\, X\in\mathcal{E}_x \,
\&\,\forall
Y\in\mathcal{E}_x,\,\alpha(Y)=\varpi_x(X,Y)\}.\end{equation}
From (\ref{eqE}), it follows that $E_x\subseteq D_x(E)$ and that we have
\begin{equation}\label{eqDE2} D_x(E)=\{(X,\alpha)+(0,\gamma)\,/\,(X,\alpha)
\in E_x,\,\gamma\in ann\,\mathcal{E}_x\}=E_x+ann\,\mathcal{E}_x.\end{equation}
Since $E_x\cap ann\,\mathcal{E}_x=ann\,\mathcal{E}'_x$,
(\ref{relptdimE}) and (\ref{eqDE2}) show that $dim\,D_x(E)=m$.
Notice that $D_x(E)\subseteq E'_x$. If $D(E)$ is differentiable and
$E$ is integrable, $\mathcal{E}$ is a generalized foliation and  $D(E)$
is integrable because of (\ref{anulpi}). \end{proof}

For instance, in the case of Example \ref{exemplufol},
$D(E)=F\oplus ann\,F$ and, in the case of Example \ref{locprod},
we get $D(E)=E\oplus T^*\mathcal{F}_2$.
\begin{example}\label{exemplucuTM}
{\rm Let $M$ be a manifold endowed with a regular, involutive,
$k$-dimensional subbundle $
\mathcal{E}\subseteq TM$ and with a $d_{tr}$-closed, truncated
$2$-form $\varpi$ of the pair $ (\mathcal{E},\mathcal{E}'=TM)$.
Then, the hypotheses of the regular case of Proposition \ref{caractnearly} are
satisfied and via (\ref{eqE}), we get the corresponding integrable, big-isotropic
structure $E_\varpi=graph\,\flat_\varpi\subseteq T^{big}M$ of rank $k$.
Let $\tilde{\mathcal{E}}$ be a complementary subbundle of $\mathcal{E}$ $(TM=
\mathcal{E}\oplus\tilde{\mathcal{E}})$
and denote by a prime and a double prime the projections of a vector
on $\mathcal{E},\tilde{\mathcal{E}}$, respectively. Then, we can extend
$\varpi$ to a $2$-form $\theta\in\Omega^2(M)$ by the formula
$$\theta(Y_1,Y_2)=\varpi(Y'_1,Y'_2)+\varpi(Y'_1,Y''_2)-\varpi(Y'_2,Y''_1)$$
and we get $E_\varpi=E_\theta$, where $E_\theta$ was defined in Example
\ref{exgraphtheta}.
We would also like to comment on the following particular case. If
$\varpi|_{\mathcal{E}\times\mathcal{E}}$ is non degenerate on each
leaf of $\mathcal{E}$ we will say that the structure
$graph\,\flat_\varpi$ is non degenerate and there exists a regular
Poisson bivector field $\Pi\in\chi^2(M)$ (Schouten-Nijenhuis
bracket $[\Pi,\Pi]=0$) with the symplectic foliation
$\mathcal{E}$. Furthermore, since the non-degeneracy of
$\varpi|_{\mathcal{E}\times\mathcal{E}}$ implies that
$ker\,\flat_\varpi=0$ and $im\,\flat_\varpi\cap
ann\,\mathcal{E}=0$ where $\flat_\varpi:\mathcal{E}\rightarrow
\mathcal{E}'^*=T^*M$, we deduce that $T^*M=im\,\flat_\varpi\oplus
ann\,\mathcal{E}$, therefore, $\varpi$ also defines a normal
bundle $Q$ of the foliation $\mathcal{E}$ which may be seen as
$Q=\mathcal{E}^{\perp_\varpi}$. Conversely, if a regular Poisson
structure $\Pi$ with symplectic foliation $\mathcal{E}$ and a
normal bundle $Q$ of $\mathcal{E}$ are given, we may extend the
leaf-wise symplectic form of $\Pi$ by the value zero on $Q$ to a
closed, truncated form $\varpi$ and we will have the
corresponding, integrable, big-isotropic structure
$graph\,\flat_{\varpi}$. Thus, a non degenerate structure
$graph\,\flat_\varpi$ is equivalent with a regular Poisson
structure together with a normal bundle of its symplectic
foliation. More exactly, if the Poisson bivector field is $P$ and
the normal bundle is $Q$, one has $E=\{(\sharp_P\lambda,\lambda)\,
/\,\lambda\in ann\,Q\}$ and $E'=E\oplus(Q\oplus Q^*$).}
\end{example}

We finish this section by indicating that, like a Dirac structure
\cite{C}, an integrable, big-isotropic structure $(E,E')$ on $M$
allows for a partial Hamiltonian formalism as follows.  A function
$f\in C^\infty(M)$ will be called a {\it Hamiltonian function} if
there exists a vector field $X_f\in\chi^1(M)$ such that $(X_f,df)\in
\Gamma E$. Similarly, if $(X_f,df)\in \Gamma E'$ $f$ is a {\it
weak-Hamiltonian function}. We will denote by $C^\infty_{Ham}(M)$
the set of Hamiltonian functions and by $C^\infty_{wHam}(M)$ the set
of weak Hamiltonian functions. The vector field $X_f$ is a {\it
Hamiltonian}, respectively {\it weak-Hamiltonian} vector field of
$f$ and any $Z\in\chi^1(M)$ which is a (weak-)Hamiltonian vector
field for some $f$ is a ({\it weak-}){\it Hamiltonian vector field}.
The fields $X^1_f,X^2_f$ are Hamiltonian vector fields of the same
function $f$ iff $X^2_f-X^1_f\in ann\,pr_{T^*M}E'$ and are
weak-Hamiltonian vector fields of $f$ iff $X^2_f-X^1_f\in
ann\,pr_{T^*M}E$. Similarly, $Z$ is Hamiltonian (weak-Hamiltonian)
for two functions $f_1,f_2$ iff $df_2-df_1\in ann\,\mathcal{E}'$
(respectively, $df_2-df_1\in ann\,\mathcal{E}$). We will denote by
$\chi_{Ham}(M),\chi_{wHam}$, respectively, the set of Hamiltonian
and weak-Hamiltonian vector fields.

Furthermore, if $f\in C^\infty_{Ham}(M)$ and $h\in
C^\infty_{wHam}(M)$ the following bracket is well defined
\begin{equation}\label{crosetfh}
\{f,h\}=X_fh=-\varpi(X_f,X_h)=-X_hf\end{equation} and does not depend
on the choice of the Hamiltonian vector fields of the functions
$f,h$. The bracket (\ref{crosetfh}) will be called the {\it Poisson
bracket} of the two functions. Formula (\ref{crosetinC}) shows that
$\{f,h\}\in C^\infty_{wHam}(M)$ and one of its weak-Hamiltonian
vector fields is $[X_f,X_h]$. Moreover, if both $f,h\in
C^\infty_{Ham}(M)$, their Poisson bracket is skew symmetric and
belongs to $C^\infty_{Ham}(M)$. Now, if we notice that
$d_{tr}\varpi(X_f,X_h,X_l)$ makes sense $\forall f,h\in
C^\infty_{Ham}(M), \forall l\in C^\infty_{wHam}(M)$ (since the
functions $\varpi(X_f,X_l)$, etc. are differentiable), the
computation done during the proof of Proposition \ref{caractnearly}
now yields $d_{tr}\varpi(X_f,X_h,X_l)=0$. This is easily seen to be
equivalent with the {\it Leibniz property}
\begin{equation}\label{Leibniz}
\{f,\{h,l\}\} =\{ \{f,h\},l\} +\{h,\{f,l\}\},
\end{equation} which restricts to the Jacobi identity on $C^\infty_{Ham}(M)$.
Therefore, $C^\infty_{Ham}(M)$ with the Poisson bracket is a Lie
algebra and $C^\infty_{wHam}(M)$ is a module over this Lie algebra.
Using the Poisson bracket  (\ref{crosetfh}) it also follows easily
that $\chi_{Ham}(M)$ is a Lie subalgebra of $\chi^1(M)$ and
$\chi_{wHam}(M)$ is a module over the former.
\section{Canonical local bases}
We will discuss local properties of a $k$-dimensional,
big-isotropic structure $E$ by constructing canonical,
local bases in the neighborhood of a fixed point
$x_0\in M$, as constructed by Dufour and Wade in the
Dirac case \cite{DW}.  In what follows the notation
is the same as in Section 2.

We begin with vectors $X^0_a,Y^0_h\in T_{x_0}M$ where $X^0_a$,
$a=1,...,dim\,\mathcal{E}_{x_0}$, is a basis of
$\mathcal{E}_{x_0}$ and $Y^0_h$, $h=1,...,dim\,\mathcal{E}'_{x_0}
-dim\,\mathcal{E}_{x_0}\stackrel{(\ref{relptdimE})}{=}m-k$, is a
basis of a complement of $\mathcal{E}_{x_0}$ in
$\mathcal{E}'_{x_0}$. Then, there exist covariant vectors
$\xi^a_0,\eta_0^h\in T^*_{x_0}M$ such that $(X^0_a,\xi_0^a)$ are
linearly independent elements of $E_{x_0}$ and
$(X^0_a,\xi_0^a),(Y^0_h,\eta_0^h)$ is a basis of a complement of
$ker(E'_{x_0}\rightarrow\mathcal{E}'_{x_0})$. Since
$ker(E'_{x_0}\rightarrow\mathcal{E}'_{x_0})=ann\,\mathcal{E}_{x_0}$,
if we add a basis $(0,\zeta^s_0)$,
$(s=1,...,dim\,ann\,\mathcal{E}_{x_0})$ of
$ann\,\mathcal{E}_{x_0}$ we get a basis of $E'_{x_0}$. Moreover,
since $ann\,\mathcal{E}_{x_0}\supseteq ann\,\mathcal{E}'_{x_0}$,
we may ask the basis $(0,\zeta^s_0)$ to consist of elements
$(0,\kappa^u_0),(0,\nu^q_0)$,
$u=1,...,dim\,ann\,\mathcal{E}'_{x_0}$,
$q=1,...,dim\,ann\,\mathcal{E}_{x_0}
-dim\,ann\,\mathcal{E}'_{x_0}\stackrel{(\ref{relptdimE})}{=}m-k$,
where $(0,\kappa^u_0)$ is a basis of $ann\,\mathcal{E}'_{x_0}$.
Then, since
$ker(E_{x_0}\rightarrow\mathcal{E}_{x_0})=ann\,\mathcal{E}'_{x_0}$,
$(X^0_a,\xi_0^a),(0,\kappa^u_0)$ is a basis of $E_{x_0}$.
Furthermore, we shall need vectors $Z_\sigma^0\in T_{x_0}M$,
$\sigma=1,...,dim\,M-dim\,\mathcal{E}'_{x_0}$
$=dim\,ann\,\mathcal{E}'_{x_0}$, which are a basis of a complement
of $\mathcal{E}'_{x_0}$ in $T_{x_0}M$.  Notice the important fact
that the indices $h,q$, on one side, and $u,\sigma$, on the other
side, have the same range.

Now we shall extend the basis
$\mathcal{B}_0=\{(X^0_a,\xi_0^a),(0,\kappa^u_0),$
$(Y^0_h,\eta_0^h),(0,\nu^q_0)\}$ to a basis of cross sections of
$E,E'$ over a neighborhood $U$ of $x_0$ in $M$; we will allow $U$ to
undergo as many restrictions as needed for the correctness of the
various constructions below without changing its name.

Clearly, we may assume that there exists a basis of $T_UM$ that
consists of local vector fields $(X_a,Y_h,Z_\sigma)$ with the
values $(X^0_a,Y^0_h,Z^0_\sigma)$ at $x_0$ and we denote by
$(\theta^a,\phi^h,\psi^\sigma)$ the corresponding, dual, local
basis of $T^*_UM$, i.e.,
$$\theta^a(X_b)=\delta^a_b,\theta^a(Y_h)=0,\theta^a(Z_\sigma)=0,$$
$$\phi^h(X_b)=0,\phi^h(Y_l)=\delta^h_l,\phi^h(Z_\sigma)=0,$$
$$\psi^\sigma(X_b)=0,\psi^\sigma(Y_k)=0,\psi^\sigma(Z_\tau)=
\delta^\sigma_\tau.$$  Accordingly, an extension  of the basis
$\mathcal{B}_0$ to a basis of $(E,E')$
over $U$ has an expression of the form
\begin{equation}\label{baza1}\begin{array}{l}
\mathcal{X}_a= (A_a^bX_b+A_a^{'h}Y_h+A_a^{''\sigma}Z_\sigma,
\alpha^a_b\theta^b+\alpha^{'a}_h\phi^h+\alpha^{''a}_\sigma\psi^\sigma),
\vspace{2mm}\\
\Xi_u= (B_u^aX_a+B_u^{'h}Y_h+B_u^{''\sigma}Z_\sigma,
\beta^u_a\theta^a+\beta^{'u}_h\phi^h+\beta^{''u}_\sigma\psi^\sigma),
\vspace{2mm}\\
\mathcal{Y}_h= (C_h^aX_a+C_h^{'l}Y_l+C_h^{''\sigma}Z_\sigma,
\gamma^h_a\theta^a+\gamma^{'h}_l\phi^l+\gamma^{''h}_\sigma\psi^\sigma),
\vspace{2mm}\\
\Theta_q= (L^a_qX_a+L_q^{'h}Y_h+L_q^{''\sigma}Z_\sigma,
\lambda^q_a\theta^a+\lambda^{'q}_h\phi^h+\lambda^{''q}_\sigma\psi^\sigma),
\end{array}\end{equation} where we use the Einstein summation convention,
$\mathcal{X}_a,\Xi_u$ is a basis
of $E|_U$, $\mathcal{Y}_h,\Theta_q$ completes the former to a
basis of $E'|_U$ and
\begin{equation}\label{condbaza1}\begin{array}{l}
A_a^b(x_0)=\delta_a^b,A_a^{'h}(x_0)=0,A_a^{''\sigma}(x_0)=0,\vspace{2mm}\\
B_u^a(x_0)=0,B_u^{'h}(x_0)=0,B_u^{''\sigma}(x_0)=0,\vspace{2mm}\\
C_h^a(x_0)=0,C_h^{'l}(x_0)=\delta_h^l,C_h^{''\sigma}(x_0)=0,\vspace{2mm}\\
L_q^a(x_0)=0,L_q^{'h}(x_0)=0,L_q^{''\sigma}(x_0)=0.\vspace{2mm}\\
\end{array}\end{equation}

We shall change this basis in order to simplify the expressions
(\ref{baza1}). However, we keep denoting the elements of the new
bases by the same letters as in (\ref{baza1}). Firstly, in view of
(\ref{condbaza1}), we may assume that the matrix $(A_a^b)$ is non
degenerate on $U$ and change the vectors $\mathcal{X}_a$ by the
matrix $(A_a^b)^{-1}$. As a result we get a basis (\ref{baza1})
where $A_a^b=\delta_a^b$. Similarly, we may get
$C^{'l}_h=\delta_h^l$. Then, the new basis may be changed by
$\Xi_u\mapsto\Xi_u-B_u^a\mathcal{X}_a$ and get a new basis
(\ref{baza1}) where $B_u^a=0$. Similarly, we may get
$C_h^a=0,L_q^a=0$, then, with the change
$\Theta_q\mapsto\Theta_q-L^{'h}_q\mathcal{Y}_h$, also get
$L^{'h}_q=0$.

Thus, any big-isotropic structure $(E,E')$ has
local bases of the form
\begin{equation}\label{baza2}\begin{array}{l}
\mathcal{X}_a= (X_a+A_a^{'h}Y_h+A_a^{''\sigma}Z_\sigma,
\alpha^a_b\theta^b+\alpha^{'a}_h\phi^h+\alpha^{''a}_\sigma\psi^\sigma),
\vspace{2mm}\\
\Xi_u= (B_u^{'h}Y_h+B_u^{''\sigma}Z_\sigma,
\beta^u_a\theta^a+\beta^{'u}_h\phi^h+\beta^{''u}_\sigma\psi^\sigma),
\vspace{2mm}\\
\mathcal{Y}_h= (Y_h+C_h^{''\sigma}Z_\sigma,
\gamma^h_a\theta^a+\gamma^{'h}_l\phi^l+\gamma^{''h}_\sigma\psi^\sigma),
\vspace{2mm}\\
\Theta_q= (L_q^{''\sigma}Z_\sigma,
\lambda^q_a\theta^a+\lambda^{'q}_h\phi^h+\lambda^{''q}_\sigma\psi^\sigma),
\end{array}\end{equation} where $(\mathcal{X}_a,\Xi_u)$ is a basis of $E$
and (\ref{condbaza1}) holds.

Furthermore, the following $g$-orthogonality conditions must be satisfied:
\begin{equation}\label{ortoinB2}
\mathcal{X}_a\perp_g\mathcal{X}_b,\, \mathcal{X}_a\perp_g\Xi_u,\,
\Xi_u\perp_g\Xi_v,\, \mathcal{X}_a\perp_g\mathcal{Y}_h,\,\end{equation}
$$\Xi_u\perp_g\mathcal{Y}_h,\,
\mathcal{X}_a\perp_g\Theta_q,\, \Xi_u\perp_g\Theta_q.$$ In particular, the
second, fifth and sixth conditions (\ref{ortoinB2}), taken at
$x_0$ give
$$\beta_a^u(x_0)=0,\;\beta_h^{'u}(x_0)=0,\;\lambda_a^q(x_0)=0.$$
This implies that the $1$-forms $\beta_\sigma^{''u}\psi^\sigma$, on
one hand, and
$\lambda^{'q}_h\phi^h+\lambda^{''q}_\sigma\psi^\sigma$, on the other
hand are linearly independent at $x_0$ and on a neighborhood $U$ of
$x_0$, which may be used for further simplifications of the basis:
i) we may linearly change $\Xi_u$ by the inverse of the matrix
$(\beta^{''u}_\sigma)$ and get a new basis (\ref{baza2}) where
$\beta^{''u}_\sigma=\delta^u_\sigma$, ii) after change i), subtract
$\sum_\sigma\alpha^{''a}_\sigma\Xi_\sigma$,
$\sum_\sigma\gamma^{''h}_\sigma\Xi_\sigma$,
$\sum_\sigma\lambda^{''q}_\sigma\Xi_\sigma$ from
$\mathcal{X}_a,\mathcal{Y}_h,\Theta_q$, respectively, and get rid of
the terms with $\psi^\sigma$ in
$\mathcal{X}_a,\mathcal{Y}_h,\Theta_q$, iii) change ii) reaches a
situation where the forms $\lambda^{'q}_h\phi^h$ are independent and
we will change the $\Theta_q$ by the inverse of the matrix
$(\lambda_h^{'q})$ and obtain $\lambda_h^{'q}=\delta_h^q$ for the
new basis, iv) subtract $\sum_l\gamma^{'h}_l\Theta_l$ from
$\mathcal{Y}_h$ and get  $\gamma^{'h}_l=0$ in the new basis, v)
since change ii) alters the coefficients of $Y_h$ in $\mathcal{Y}_h$
and adds terms in $Y_h$ to $\Theta_q$, we correct that by changing
the new $\mathcal{Y}_h$ with the corresponding coefficient matrix
(which is non degenerate on $U$) and by subtracting the necessary
linear combination of $\mathcal{Y}_h$ from $\Theta_q$. The result is
a local basis of $(E,E')$ that looks as follows
\begin{equation}\label{baza2'}\begin{array}{l}
\mathcal{X}_a= (X_a+A_a^{'h}Y_h+A_a^{''\sigma}Z_\sigma,
\alpha^a_b\theta^b+\alpha^{'a}_h\phi^h),
\vspace{2mm}\\
\Xi_u= (B_u^{'h}Y_h+B_u^{''\sigma}Z_\sigma,
\beta^u_a\theta^a+\beta^{'u}_h\phi^h+\psi^u),
\vspace{2mm}\\
\mathcal{Y}_h= (Y_h+C_h^{''\sigma}Z_\sigma,
\gamma^h_a\theta^a),
\vspace{2mm}\\
\Theta_q= (L_q^{''\sigma}Z_\sigma,
\lambda^q_a\theta^a+\phi^q),
\end{array}\end{equation} where (\ref{condbaza1}) are still valid.

It is easy to see that, if the vector fields $X_a,Y_h,Z_\sigma$
are fixed, there is only one basis of $(E,E')$ which is of the
form (\ref{baza2'}). For this reason we will say that the basis
(\ref{baza2'}) is a {\it canonical, local basis} of the
big-isotropic structure $E$.

Now, we consider the integrable case. Let $\mathcal{U}$ be a
neighborhood of the point $x_0$ on the characteristic leaf
$\mathcal{S}$ through $x_0$. Then
$dim\,\mathcal{E}|_{\mathcal{U}}=const.$ and, in view of
(\ref{relptdimE}), $dim\,\mathcal{E}'|_{\mathcal{U}}=const.$ too.
Furthermore, on the neighborhood $U$ of $x_0$ in $M$ there are
coordinates $(x^a,y^h,z^\sigma)$ such that $x_0$ has the
coordinates $(0,0,0)$, the equations of $\mathcal{U}$ are
$y^h=0,z^\sigma=0$,
$\mathcal{E}|_{\mathcal{U}}=span\{(\partial/\partial
x^a)|_{\mathcal{U}}\}$ and
$\mathcal{E}'|_{\mathcal{U}}=span\{(\partial/\partial
x^a)|_{\mathcal{U}}, (\partial/\partial y^h)|_{\mathcal{U}}\}$.
Indeed, we may assume that $U$ is a tubular neighborhood of
$\mathcal{U}$ where the tangent space of the tubular fibers at the
points of $\mathcal{U}$ is a direct sum of a complementary space
of $\mathcal{E}$ in $\mathcal{E}'$ and a complementary space of
$\mathcal{E}'$ in $TM$; then, take $x^a$ coordinates on
$\mathcal{U}$ and $y^h,z^\sigma$ coordinates along the tubular
fibers such that $(\partial/\partial y^h)|_{\mathcal{U}}$ span the
chosen complement of $\mathcal{E}$ in $\mathcal{E}'$ and
$(\partial/\partial z^\sigma)|_{\mathcal{U}}$ span the chosen
further complement in $T_{\mathcal{U}}M$. The tubular structure of
the neighborhood  $U$ will be important and we will denote by
$\mathcal{F}$ the foliation of $U$ by the tubular fibers, for
later use. It is easy to see that we can construct bases
(\ref{baza1}) where
\begin{equation}\label{auxbaza3} X_a=\frac{\partial}{\partial x^a},\,
Y_h=\frac{\partial}{\partial y^h}+\chi^\sigma_h
\frac{\partial}{\partial
z^\sigma},\,Z_\sigma=\frac{\partial}{\partial
z^\sigma}\;\;(\chi^\sigma_h|_{\mathcal{U}}=0).\end{equation} If these values
are inserted in (\ref{baza2'}), the result takes the following form
(with new coefficients):
\begin{equation}\label{baza3}\begin{array}{l}
\mathcal{X}_a= (\frac{\partial}{\partial
x^a}+A_a^{'h}\frac{\partial}{\partial
y^h}+A_a^{''\sigma}\frac{\partial}{\partial z^\sigma},
\alpha^a_bdx^b+\alpha^{'a}_hdy^h),
\vspace{2mm}\\
\Xi_u= (B_u^{'h}\frac{\partial}{\partial
y^h}+B_u^{''\sigma}\frac{\partial}{\partial z^\sigma},
\beta^u_adx^a+\beta^{'u}_hdy^h+ dz^u),
\vspace{2mm}\\
\mathcal{Y}_h= (\frac{\partial}{\partial
y^h}+C_h^{''\sigma}\frac{\partial}{\partial z^\sigma},
\gamma^h_adx^a),
\vspace{2mm}\\
\Theta_q= (L_q^{''\sigma}\frac{\partial}{\partial z^\sigma},
\lambda^q_adx^a+dy^q),
\end{array}\end{equation} where (\ref{condbaza1}) holds
$\forall x\in\mathcal{U}$ and (\ref{ortoinB2})  holds on $U$,
which means that we have
\begin{equation}\label{ortoconcret} \begin{array}{l} \alpha^{'a}_h+\gamma^h_a=0,\,
\lambda^q_a+A^{'q}_a=0,\,\beta^{'u}_h+C^{''u}_h=0,\,L^{''u}_q+B^{'q}_u=0,
\vspace{2mm}\\ \beta_a^u+A^{''u}_a+\alpha^{'a}_hB^{'h}_u+\beta^{'u}_hA^{'h}_a=0,\,
B^{''u}_v+B^{''v}_u+\beta^{'u}_hB^{'h}_v+\beta^{'v}_hB^{'h}_u=0,\vspace{2mm}\\
\alpha^a_b+\alpha^b_a+\alpha^{'a}_hA^{'h}_b+\alpha^{'b}_hA^{'h}_a=0.
\end{array}\end{equation}

We also notice that, if we are interested in bases of $E'$,
without requiring them to be a prolongation of a basis of $E$, we
may repeat the subtraction trick between $\mathcal{X}_a,\Xi_u$ and
$\Theta_q$ and between $\mathcal{X}_a,\Xi_u$ and $\mathcal{Y}_h$
as well and get a basis of the form (\ref{baza3}) with the
supplementary conditions \begin{equation}\label{baza4}
\alpha^{'a}_h=0,\;\beta^{'u}_h=0,\;A^{'h}_a=0,\;B^{'h}_u=0;\end{equation}
the new pairs  $\mathcal{X}_a,\Xi_u$ may not belong to $E$ any
more and (\ref{ortoconcret}) does not hold.

It follows easily that, if the local coordinates $(x^a,y^h,z^\sigma)$ such that
\begin{equation}\label{coordcanonice} y^h|_{\mathcal{U}}=0,
\,z^\sigma|_{\mathcal{U}}=0, \, \left.\frac{\partial}{\partial y^h}
\right|_{\mathcal{U}}\in(\mathcal{E}'\backslash\mathcal{E})|_{\mathcal{U}},\,
\left.\frac{\partial}{\partial z^\sigma}\right|_{\mathcal{U}}
\in(TM\backslash\mathcal{E}')|_{\mathcal{U}}\end{equation}
are chosen, the basis of the form (\ref{baza3}) where
$\mathcal{X}_a,\Xi_u\in E$ and $\mathcal{Y}_h,\Theta_q\in
E'\backslash E$ is unique. A similar fact holds for the basis of
$E'$ which satisfies the conditions (\ref{baza4}). Accordingly,
the basis (\ref{baza3}) will be called a {\it canonical, local
basis} of the integrable, big-isotropic structure $(E,E')$ and the
basis where (\ref{baza4}) also holds is a {\it canonical basis} of
$E'$.
\begin{rem} {\rm In the Dirac case $k=m$, the basis (\ref{baza3})
reduces to
\begin{equation}\label{canonicDW}
\mathcal{X}_a= (\frac{\partial}{\partial
x^a}+A_a^{''\sigma}\frac{\partial}{\partial z^\sigma},
\alpha^a_bdx^b), \;
\Xi_u= (B_u^{''\sigma}\frac{\partial}{\partial z^\sigma},
\beta^u_adx^a+ dz^u),\end{equation}
which are the formulas given in \cite{DW}. In the almost Dirac
case, we have the corresponding formulas deduced from
(\ref{baza2'})
\begin{equation}\label{aproapeDW}
\mathcal{X}_a= (X_a+A_a^{''\sigma}Z_\sigma,
\alpha^a_b\theta^b),\;
\Xi_u= (B_u^{''\sigma}Z_\sigma,
\beta^u_a\theta^a+\psi^u).\end{equation} For (\ref{canonicDW}) and
(\ref{aproapeDW}) the conditions (\ref{ortoconcret}) reduce to
\begin{equation}\label{ortoinDW} \alpha_a^b+\alpha_b^a=0,\,
B^{''\sigma}_u+B_{\sigma}^{''u}=0,\,\beta_a^{u}+A^{''u}_a=0.
\end{equation} As an example, we use these formulas for a straightforward proof
of the fact that, in the two-dimensional case $m=2$, any almost
Dirac structure $D$ is Dirac. Indeed, the points of $M$ may be
classified into three classes $(0,1,2)$ where
$dim(pr_{TM}D)=0,1,2$, respectively. Obviously, any point of class
$2$ has a neighborhood $U$ of class $2$ and $D$ is integrable on
$U$. Formulas (\ref{aproapeDW}) and (\ref{ortoinDW}) show that any
point of class $1$ is regular, hence, it also has a neighborhood
where $D$ is integrable. If a point of class $0$ has a
neighborhood $U$ of points of class $0$, $D$ is integrable on $U$.
Finally, if a point $x_0$ of class $0$ has no such neighborhood,
every neighborhood of $x_0$ has points that are of class $2$
(necessarily), hence, $x_0$ is the limit of points that have
neighborhoods where $D$ is closed by Courant brackets and, by
continuity, $x_0$ also has a neighborhood where $D$ is integrable.
A similar analysis, where it is simpler to use Proposition
\ref{caractnearly} instead of the local, canonical bases shows that a
big-isotropic structure $E$ with $k=1,m=2$ must be integrable.}
\end{rem}
\section{Local geometry of a big-isotropic structure}
We begin with a preparatory discussion of the operation of pullback
of a big-isotropic structure by a mapping (see \cite{{BR},{C}} for
the Dirac case), which is of a more general interest. Let
$f:N^n\rightarrow M^m$ be a mapping of manifolds, $x\in N, y=f(x)$ a
pair of corresponding points, and $E$ an arbitrary vector subbundle
of $T^{big}M$. We denote by $f_*$ the differential of $f$ and by
$f^*$ the transposed mapping of $f_*$. Then
\begin{equation}\label{pullback} f^*(E_y) = \{(X,f^*\alpha)\,/\,
X\in T_xN,\,\alpha\in T_y^*M,\end{equation} $$(f_*X,\alpha)\in E_y\}$$
is the pullback of $E$ at $x$.
\begin{prop}\label{proppull} Let $E$ be a rank $k$, big-isotropic subbundle of
$T^{big}M$. Then $f^*(E_y)$ is isotropic in $T^{big}_xN$ and its
$g_N$-orthogonal space is $(f^*(E_y))'=f^*(E'_y)$. Furthermore, if
$f$ is an embedding of $N$ in $M$, if $E$ is an integrable,
big-isotropic structure on $M$ and if $f^*E=\cup_{y\in M}f^*(E_y)$
is a differentiable subbundle of $T^{big}N$ then $f^*E$ is an
integrable, big-isotropic structure on $N$. Finally, the condition
\begin{equation}\label{condpb} E_y\cap ker\,f^*_y = E'_y\cap
ker\,f^*_y
\end{equation} characterizes the situation where the dimension of
$f^*(E_y)$ is $n-m+k$. \end{prop}
\begin{proof} It is easy to check that $f^*(E_y)$ is an isotropic
subspace of $(T_x^{big}N,g_N)$ and we will compute its dimension.
Obviously, we have
\begin{equation}\label{eqauxptf} dim\,ker\,f_{*x}=n-rank_xf,\,
dim\,ker\,f^{*}_y=m-rank_xf.\end{equation} Then, define the
space
\begin{equation}\label{eqS} S_y= \{(f_*X,\alpha)\in E_y\,/\, X\in
T_xN,\,\alpha\in T^*_yM\} = E_y\cap(im\,f_*\oplus T^*_yM),
\end{equation} and notice that the correspondence $(f_*X,\alpha)
\mapsto (X,f^*\alpha)$ produces an isomorphism
\begin{equation}\label{isomorphaux} S_y/(S_y\cap ker\,f^*_y)\approx
f^*(E_y)/ker\,f_{*x}.\end{equation} Therefore,
\begin{equation}\label{dimaux1} dim\,f^*(E_y)=dim\,ker\,f_{*x}
+dim\,S_y- dim(S_y\cap ker\,f^*_y)\end{equation}
$$=dim\,ker\,f_{*x}+dim\,S_y- dim(E_y\cap ker\,f^*_y)$$ (the equality
$S_y\cap ker\,f^*_y=E_y\cap ker\,f^*_y$ follows from the definition
of $S_y$ since $ker\,f^*_y\subseteq im\,f_{*x}\oplus T^*_yM$).

Furthermore, we notice the equality
\begin{equation}\label{dimaux2} im\,f_{*x}\oplus
T^*_yM=(ker\,f^*_y)^{\perp_{g_M}},\end{equation} which follows since
the left hand side is included in the right hand side and the two
spaces have the same dimension by (\ref{eqauxptf}). Accordingly, we
have
\begin{equation}\label{auxS} S_y=(E'_y)^{\perp_{g_M}}\cap
(ker\,f^*_y)^{\perp_{g_M}}=(E'_y+ker\,f^*_y)^{\perp_{g_M}},\end{equation}
whence, using the classical relation between the dimensions of the sum and
intersection of two linear subspaces, we get
\begin{equation}\label{auxdimS} dim\,S_y= 2m-dim(E'_y+
ker\,f^*_y)\end{equation} $$=dim\,E_y - dim\,ker\,f^*_y + dim(E'_y\cap ker\,f^*_y).$$

If we combine (\ref{dimaux1}), (\ref{auxdimS}) and
(\ref{eqauxptf}) we get
\begin{equation}\label{auxfinal} dim\,f^*(E_y) = n-m+k +
dim(E'\cap ker\,f^*_y)-dim(E\cap ker\,f^*_y).\end{equation}

The same procedure with the roles of $E$ and $E'$ interchanged uses
the space $S'_y=E'_y\cap(im\,f_*\oplus T^*_yM)$, which (like in
(\ref{auxdimS})) has the dimension
\begin{equation}\label{dimS'} dim\,S'_y=dim\,E'_y- dim\,ker\,f^*_y
+dim(E_y\cap ker\,f^*_y),\end{equation}
and leads to
\begin{equation}\label{auxfinal2} dim\,f^*(E'_y) = n+m-k -
dim(E'\cap ker\,f^*_y)+dim(E\cap ker\,f^*_y).\end{equation} Hence,
$dim\,f^*(E_y)+dim\,f^*(E'_y)=2n$ and, since it is easy to check
the $g_N$-orthogonality of these two spaces, we have
$f^*(E'_y)=(f^*(E_y))^{\perp_{g_N}}$.

For the second assertion of the proposition, it suffices to notice
that $X\in T_xN$ belongs to $pr_{T_xN}f^*(E_y)$ iff $f_*X\in
pr_{T_yM}E_y$. If $f$ is an embedding and if $f^*E$ is a
differentiable subbundle (in particular, $dim\,f^*(E_y)=const.$), a
field $X\in\Gamma (pr_{T_xN}(f^*E))$ has an $f$-related field
$f_*X\in\Gamma (pr_{T_yM}E)$. Moreover, since $dim\,ker\,f^*_y=
m-n=const.$, a cross section of $f^*E$ must be of the form
$(X,f^*\alpha)$ where $(f_*X,\alpha)$ is a differentiable cross
section of $E$, and the same holds for $f^*E'$ and $E'$. If these
facts are taken into consideration, a straightforward examination of
the Courant brackets shows that the integrability conditions for $E$
imply the integrability conditions for $f^*E$.

Finally, hypothesis (\ref{condpb}) is equivalent with $dim(E'\cap
ker\,f^*_y)=dim(E\cap ker\,f^*_y)$ and, by (\ref{auxfinal}), we have
the required dimension for $f^*(E_y)$ iff (\ref{condpb})
holds.\end{proof}
\begin{rem}\label{obscodim} {\rm1) By (\ref{auxfinal}), $dim\,f^*(E_y)=const.$
iff $$dim(E'_y\cap ker\,f^*_y) - dim(E_y\cap ker\,f^*_y)
=const.$$ By (\ref{auxdimS}) and (\ref{dimS'}), this condition is equivalent with
$$dim\,S'_y-dim\,S_y=2(m-k)-(dim(E'\cap
ker\,f^*_y)+dim(E\cap ker\,f^*_y))=const.$$

2) The algebraic content of Proposition \ref{proppull} holds for
any linear mapping $l:V\rightarrow W$ between linear spaces and
any isotropic subspace $E\subseteq W\oplus W^*$. The significance
of the equality $dim\,f^*(E_y)=n-m+k$ is the codimensional
invariance property $n-dim\,f^*(E_y)=m-dim\,E_y$. Hypothesis
(\ref{condpb}) holds in the Dirac case $(E=E')$ and in the case
where $f$ is a submersion $(ker\,f^*_y=0)$.}\end{rem}
\begin{corol}\label{cazdedif} With the notation of Proposition
{\rm\ref{proppull}} and of its proof, if $f$ is an embedding and if
$dim\,S_y=const.$, $dim\,S'_y=const.$ the pullbacks $f^*E,f^*E'$ are
differentiable.
\end{corol}
\begin{proof} By the definition of $S_y,S'_y$ (see (\ref{eqS})), if the
dimensions of these spaces do not depend on $x$, $S_y,S'_y$  are
differentiable with respect to $x$; then, by Remark \ref{obscodim} 1)
$dim\,f^*(E_y)=const.$ and $dim\,f^*(E'_y)=const.$ On the other hand,
since $f$ is an embedding, $ker\,f_{*x}=0$ and formula
(\ref{isomorphaux}) yields
\begin{equation}\label{eqscuf} S_y/(E_y\cap ker\,f^*_y)\approx
f^*(E_y), S'_y/(E'_y\cap ker\,f^*_y)\approx f^*(E'_y).\end{equation}
Formula (\ref{eqscuf}) and the constant dimensions of the spaces
therein imply the differentiability of $f^*E,f^*E'$ and Proposition
\ref{proppull} allows us to conclude. \end{proof}

Now, we shall discuss some local properties of an integrable,
big-isotropic structure $E$.
\begin{prop}\label{induspefoaie} The pullback of an integrable,
big-isotropic structure $E$ to a characteristic leaf $\mathcal{S}$
is the same as the pullback of its Dirac extension $D(E)$ and it
is a presymplectic structure of $\mathcal{S}$.\end{prop}
\begin{proof} Consider the neighborhoods $U,\mathcal{U}$
where one has the canonical basis (\ref{baza3}) and
$\mathcal{S}\cap U$ has the equations $y^h=0,z^\sigma=0$. Then, at
the points of $\mathcal{U}$, $\gamma\in\Omega^1(M)$ belongs to
$ann\mathcal{E}|_{\mathcal{U}}$ iff
$\gamma|_{\mathcal{U}}=\gamma_hdy^h+\gamma_\sigma dz^\sigma$. This
implies that $i^*\gamma=0$ $(i:\mathcal{U}\rightarrow M)$, whence,
by (\ref{eqDE2}), $i^* D(E)=i^*E$. Now, (\ref{condbaza1}) and
(\ref{ortoconcret}), which hold for the coordinate values
$(x^a,0,0)$, give $\beta^u_a(x^c,0,0)=0, \alpha_a^b(x^c,0,0)=
-\alpha^a_b(x^c,0,0)$ and, accordingly, (\ref{pullback}) shows
that $$i^*E=\{(f_a\frac{\partial}{\partial
x^a},f_a\alpha^a_b(x^c,0,0)dx^b)\,/\, f_a=f_a(x^c)\}.$$ Therefore,
$i^*E$ is the presymplectic structure defined by the $2$-form of
components $\alpha_b^a(x^c,0,0)$, which is just
$\varpi|_{\mathcal{E}\times\mathcal{E}}$.
\end{proof}
\begin{rem}\label{obsDE} {\rm The canonical basis (\ref{baza3}) also provides a
local expression of the Dirac extension $D(E)$ over the neighborhood $U$.
Indeed, (\ref{baza3}) shows that the annihilator of $span\{pr_{TM}\mathcal{X}_a\}$
is spanned by the $1$-forms $dy^h-A^{'h}_adx^a,dz^\sigma-A^{''\sigma}_adx^a$
and $ann\,\mathcal{E}$ consists of the forms $$\gamma=\varphi_h(dy^h-A^{'h}_adx^a) +
\psi_\sigma(dz^\sigma-A^{''\sigma}_adx^a)$$ where
\begin{equation}\label{auxfipsi} \varphi_hB^{'h}_u+\psi_\sigma B^{''\sigma}_u=0.
\end{equation} Therefore, \begin{equation}\label{DEU}
D_U(E)=\{f^a\mathcal{X}_a+s^u\Xi_u+(0,\varphi_h (dy^h-A^{'h}_adx^a)  \end{equation} $$+
\psi_\sigma(dz^\sigma-A^{''\sigma}_adx^a))\}$$ where $f^a,s^u,\varphi_h,\psi_\sigma
\in C^\infty(U)$ and (\ref{auxfipsi} ) holds. In particular, if the structure
$E$ is regular then $B^{'h}_u\equiv0,
B^{''\sigma}_u\equiv0$ and $D(E)$ is a differentiable, Dirac structure.}\end{rem}

Now, at $x_0$ , we consider the local transversal submanifold
$\mathcal{Q}_0$ of the characteristic leaf $\mathcal{S}$ given
by the equations $x^a=0$ and denote by
$\iota:\mathcal{Q}_0\rightarrow M$ the corresponding embedding.
For this embedding, (\ref{baza3}) shows that
the spaces $S,S'$ of Corollary \ref{cazdedif} are given by
$$S=span\{\Xi_u|_{\mathcal{Q}_0}\},\;\;S'=span\{\Xi_u|_{\mathcal{Q}_0},
\mathcal{Y}_h|_{\mathcal{Q}_0},\Theta_q|_{\mathcal{Q}_0}\},$$
hence, $S,S'$ have a constant dimension. Then, Corollary \ref{cazdedif}
tells us that $\mathcal{Q}_0$ has an
induced, integrable, big-isotropic structure
$E^{tr}_{x_0}=\iota^*E$, which will be called the {\it transversal
structure} of $E$ at $x_0$.

Furthermore, with (\ref{eqS} ) and (\ref{baza3}), it follows that
$E\cap ker\,\iota^*= S\cap ker\,\iota^*=0$ and we see that $E_x^{tr}$ is
isomorphic with the bundle $S|_{\mathcal{Q}_0}$. In particular,
$\Xi_u$ $({\rm mod}\,x^a=0)$ is a basis of the transversal structure
$E^{tr}_{x_0}$ along $\mathcal{Q}_0$. We also notice that
$E^{tr}\cap (T\mathcal{F}\oplus ann\,T\mathcal{F})=0$ and, in particular,
$E^{tr}_{x_0}$ is of the graph type (it is the graph of a mapping
$\Lambda\rightarrow T{\mathcal Q}_0$, where $\Lambda$ is a field of
subspaces of $T^*\mathcal{Q}_0$).

Any local transversal submanifold $\mathcal{T}_0$ of $\mathcal{U}$
at $x_0$ inherits a transversal structure since there exists a
tubular neighborhood such that $\mathcal{Q}_0=\mathcal{T}_0$.
Hence, we may speak of the transversal structure in a generic way,
which, however, does not mean that the transversal structures
defined on different submanifolds $\mathcal{T}_0$ are equivalent.
We shall address this question but we need more preparations
first. The following considerations are inspired by the case of
foliation-coupling Dirac structures \cite{VJGP} and may be used in
a discussion of the coupling between a big-isotropic structure and
a foliation, which we do not intend to develop.

Recall that $U$ has the tubular foliation $\mathcal{F}$ of
equations $x^a=const.$. Define the field of subspaces
\begin{equation}\label{pseudoE} H(E,\mathcal{F})=
\{Z\in TM\,/\,\exists\alpha\in ann\,T\mathcal{F},\,(Z,\alpha)\in E\}.
\end{equation} Using the basis (\ref{baza3} we get
\begin{equation}\label{pseudoE2}H(E,\mathcal{F})=
\{f^apr_{TM}\mathcal{X}_a\,/\,\sum_af^a\alpha^{'a}_u=0,\,f^a\in C^\infty(U)\},
\end{equation} therefore, $H(E,\mathcal{F})\cap T\mathcal{F}=0$.
Generally, $H(E,\mathcal{F})$ may not have a constant dimension;
we will say that $H(E,\mathcal{F})$ is the {\it pseudo-normal}
bundle of $\mathcal{F}$. $H(E,\mathcal{F})$ is a true normal
bundle, i.e.,
\begin{equation}\label{HnormalF}
T_UM=H(E,\mathcal{F})\oplus T\mathcal{F},\end{equation} iff
\begin{equation}\label{eqcoupling} \alpha^{'a}_h=0 \hspace{5mm}
(\stackrel{(\ref{ortoconcret})}{\Leftrightarrow}\,\gamma_a^h=0).\end{equation}

We would like to notice that condition (\ref{eqcoupling}) has other
interesting interpretations too.
One of them is obtained if we define the field of subspaces
\begin{equation}\label{subcampdual} \mathcal{H}(E',\mathcal{F})=\{\theta\in T^*M\,/\,
\exists Z\in T\mathcal{F}, (Z,\theta)\in E'\},\end{equation} which
we call the {\it pseudo-conormal bundle of $\mathcal{F}$ modulo}
$E$. Using (\ref{baza3}), it follows that
$\theta\in\mathcal{H}(E',\mathcal{F})$ iff
$$\theta=\varphi_u(\beta^u_adx^a+\beta^{'u}_hdy^h+ dz^u) + \psi_h\gamma^h_adx^a
+\mu_q(\lambda^q_adx^a+dy^q).$$
Accordingly, we see that
\begin{equation}\label{HrondsiannF} T^*_UM=\mathcal{H}(E',\mathcal{F})+
ann\,T\mathcal{F},\; \mathcal{H}(E',\mathcal{F})\cap ann\,T\mathcal{F}=
span\{\gamma_a^hdx^a\}.\end{equation}
Therefore, we have
\begin{equation}\label{coupling} T^*_UM=\mathcal{H}(E',\mathcal{F})
\oplus ann\,T\mathcal{F}\end{equation} iff (\ref{eqcoupling}) holds.

For another interpretation of condition (\ref{eqcoupling}) let us
define the field of subspaces
\begin{equation}\label{pseudonormal} H(E',\mathcal{F})=
\{Z\in TM\,/\,\exists\alpha\in ann\,T\mathcal{F},\,(Z,\alpha)\in E'\}.
\end{equation} Notice that $H(E,\mathcal{F})
\subseteq H(E',\mathcal{F})$.
From (\ref{baza3}), we get \begin{equation}\label{eqHE'}
H(E',\mathcal{F})=span\{pr_{TM}\mathcal{X}_a-
\sum_{q,\sigma}\alpha^{'a}_qL^{''\sigma}_q\frac{\partial}{\partial z^\sigma},
pr_{TM}\mathcal{Y}_h\},\end{equation} Therefore,
$\mathcal{E}'|_U=H(E',\mathcal{F})+T\mathcal{F}$ and
$H(E',\mathcal{F})\cap T\mathcal{F}$ is the (trivial)
$(m-k)$-dimensional bundle with the basis
$\{pr_{TM}\mathcal{Y}_h\}$. Together with (\ref{eqHE'}), this
implies
\begin{equation}\label{dulptHE'} (ann\,\mathcal{F})^*\approx H(E',\mathcal{F})
/(H(E',\mathcal{F})\cap T\mathcal{F}).\end{equation}

On the other hand, let us recall the truncated $2$-form
$\varpi:\mathcal{E}\times\mathcal{E}' \rightarrow \mathds{R}$ and
consider its {\it second flat morphism} $\flat'_{\varpi}:\mathcal{E}'
\rightarrow\mathcal{E}^*$ given by $(\flat'_{\varpi}Y)(X)=\varpi(X,Y)$.
It follows that condition (\ref{eqcoupling}) holds iff
\begin{equation}\label{cupl2} H(E',\mathcal{F})\cap T\mathcal{F}\subseteq ker\,
\flat'_\varpi. \end{equation} This again is an interpretation of (\ref{eqcoupling}).
Notice also that if (\ref{eqcoupling}) holds then $\forall X\in
H(E',\mathcal{F})$ the corresponding $1$-form $\alpha$ such that
$(X,\alpha)\in E'$ is uniquely defined.

The interesting fact that follows from (\ref{HnormalF}) is that
$E|_U$ decomposes along $\mathcal{F}$ and $H(E,\mathcal{F})$.
Indeed, using the canonical basis (\ref{baza3}) we see that
(\ref{HnormalF}) implies
$$T^*\mathcal{F}\approx ann\,H(E,\mathcal{F})
=span\{dy^h-A^{'h}_adx^a,\,dz^\sigma-A^{''\sigma}_adx^a\}$$ and
$$E\cap(T\mathcal{F}\oplus T^*\mathcal{F})=span\{\Xi_u\}.$$ Thus,
this intersection produces the transversal structures on the
fibers of $\mathcal{F}$.  We also have $H^*(E,\mathcal{F})\approx
ann\,T\mathcal{F}$ and
$$E\cap(H(E,\mathcal{F})\oplus H^*(E,\mathcal{F}))=span\{\mathcal{X}_a\}.$$
Therefore, the following decomposition holds
\begin{equation}\label{descluiE} E|_U=[E\cap(T\mathcal{F}\oplus T^*\mathcal{F})]
\oplus[E\cap(H(E,\mathcal{F})\oplus H^*(E,\mathcal{F}))].\end{equation}

This property justifies the following definition.
\begin{defin}\label{localdec} {\rm An integrable,
big-isotropic structure $E$ on $M$ is
called {\it locally decomposable} if each point $x\in M$ has a
tubular neighborhood $U$ of the characteristic slice (a
neighborhood of the characteristic leaf) $\mathcal{U}$ through
$x$, with the tubular foliation $\mathcal{F}$, where the
decomposition (\ref{HnormalF}) holds. If (\ref{HnormalF}) holds
for any tubular neighborhood, $E$ will be called {\it strongly,
locally decomposable}.}\end{defin}

All the Dirac structures are strongly, locally decomposable \cite{DW};
a non-Dirac example follows. \begin{example}\label{exstrong}
{\rm Take $M=\mathds{R}^3$ with coordinates $(x,y,z)$ and
$$E=span\{(\frac{\partial}{\partial x}, 0),(0,dz)\}.$$ $E$ is a
regular, integrable, big-isotropic structure of dimension $2$ with
$$E'=E\oplus span\{(\frac{\partial}{\partial y},0),(0,dy)\}.$$
The given basis is canonical and local decomposability holds.
New tubular coordinates are defined by
$$\tilde x=\tilde x(x,y,z),\,\tilde y=\tilde y(x,y,z),\,\tilde z=\tilde z(x,y,z),$$
where (see (\ref{coordcanonice})  $$\tilde y(x,0,0)=0,\,\tilde z(x,0,0)=0,\,
\left.\frac{\partial\tilde y}{\partial x}\right|_{(x,0,0)}=0,\, \left.
\frac{\partial\tilde z}{\partial x}\right|_{(x,0,0)}=0,\,$$
$$\left.\frac{\partial\tilde z}{\partial y}\right|_{(x,0,0)}=0,\,
\frac{\partial\tilde x}{\partial x}\neq0,\,
\frac{\partial\tilde y}{\partial y}\neq0,\,
\frac{\partial\tilde z}{\partial z}\neq0.$$
If we make this change in the generators of $E$ and produce a
canonical basis out of the result we see that (\ref{eqcoupling}),
hence local decomposability, also holds for the new tubular
neighborhood.}\end{example}

The following example shows that local decomposability of an
integrable, big-isotropic structure does not imply strong, local
decomposability.
\begin{example}\label{exdenedecomp} {\rm Take $M=\mathds{R}^5$ with the canonical
coordinates $(x^1,x^2,y^1,y^2,z)$. Define
\begin{equation}\label{eq1inexn} E= span\{\mathcal{X}_1=
(\frac{\partial}{\partial x^1}, dx^2+dy^1),\end{equation}
$$\mathcal{X}_2=(\frac{\partial}{\partial x^2},- dx^1+
dy^2),\,\Xi_1=(0,dz)\}.$$
This is a regular, $3$-dimensional, integrable, big-isotropic structure
and the orthogonal bundle $E'$ has the supplementary generators
\begin{equation}\label{eq2inexn} \mathcal{Y}_1=
(\frac{\partial}{\partial y^1},- dx^1),\, \mathcal{Y}_2=
(\frac{\partial}{\partial y^2},- dx^2),\,\Theta_1=(0,dy^1),\,
\Theta_2=(0,dy^2).\end{equation}
The characteristic leaf through the origin is  the $(x^1,x^2)$-plane,
$\mathds{R}^5$ may be seen as a tubular neighborhood of this leaf and the basis given by (\ref{eq1inexn}), (\ref{eq2inexn}) is canonical. Accordingly, the local decomposability property does not hold for the tubular foliation that consists of the family of $3$-dimensional planes with coordinates $(y^1,y^2,z)$.
Now, consider the following coordinate transformation
$$\tilde x^1=x^1-y^2,\, \tilde x^2=x^2+y^1,\,
\tilde y^1=y^1,\,\tilde y^2=y^2,\,\tilde z=z.$$ If we express the pairs
(\ref{eq1inexn}), (\ref{eq2inexn}) and then produce a canonical basis
of $(E,E')$ with respect to the new coordinates we get
$$\tilde{\mathcal{X}}_1=(\frac{\partial}{\partial\tilde x^1}, d\tilde x^2),\,
\tilde{\mathcal{X}}_2=(\frac{\partial}{\partial\tilde x^2},- d\tilde x^1),\,
\tilde\Xi=(0,d\tilde z),$$
$$\tilde{\mathcal{Y}}_1=(\frac{\partial}{\partial\tilde y^1},0),\,
\tilde{\mathcal{Y}}_2=(\frac{\partial}{\partial\tilde y^2},0),\,
\tilde\Theta_1=(0,d\tilde y^1),\,\tilde\Theta_2=(0,d\tilde y^2).$$
Thus, with respect to the tubular fibers defined by the new coordinates
the local decomposability property holds.}\end{example}

It is known that the transversal structure of a Dirac structure is
well defined up to a natural equivalence \cite{DW}.  We will show
that the same holds for the integrable, big-isotropic structures
that are strongly, locally decomposable. The basis (\ref{baza3})
also yields a transversal structure on each submanifold
$\mathcal{Q}_{x}$ defined by $x^a=x^a(x)$ $(x\in\mathcal{U})$.  We
will denote by $E^{tr}$ the family of transversal structures
$E^{tr}_x$, $x\in\mathcal{U}$ and prove the following lemma.
\begin{lemma}\label{lemaptisotopie} Let $E$ be a locally decomposable,
integrable, big-isotropic structure on $M$, $x_0$ a point of $M$ and $U$ a
tubular neighborhood of $x_0$ with the tubular foliation $\mathcal{F}$.
Then, any $\mathcal{F}$-projectable vector field $Z\in H(E,\mathcal{F})$
is an infinitesimal automorphism of the transversal structure $E^{tr}$.\end{lemma}
\begin{proof} Since $Z$ is $\mathcal{F}$-projectable, the flow of $Z$ sends
leaves of $\mathcal{F}$ to leaves of $\mathcal{F}$. As explained
earlier, the transversal structure $E^{tr}$ on the leaves of
$\mathcal{F}$ is induced by the vector bundle $S=span\{\Xi_u\}$
(see (\ref{baza3})).  Equivalently, $E^{tr}$ is induced also by
the bundle $\tilde{E}^{tr}=S\oplus ann\,T\mathcal{F}$ (obviously,
the intersection of the terms is zero and $\tilde{E}^{tr}$ is an
integrable, big-isotropic structure on $U$ of the same dimension
$k$ like $E$). Thus, the conclusion will be obtained if we show
that $Z$ is an infinitesimal automorphism of $\tilde{E}^{tr}$. If
we denote
\begin{equation}\label{brievety}\mathcal{X}_a=(V_a,\nu^a),\,\Xi_u=(W_u,\xi^u),
\end{equation}
the projectability of $Z$ is equivalent with asking
$Z=\zeta^a(x^b)V_a$ and $Z$ is an infinitesimal automorphism of
$\tilde{E}^{tr}$ iff
\begin{equation}\label{auxZ} (L_{V_a}W_u,L_{V_a}\xi^u)\in\Gamma\tilde{E}^{tr}.
\end{equation}

The integrability of $E$ implies the existence of local functions
$f^a,\varphi^u$ such that
$$[\mathcal{X}_a,\Xi_u] = ([V_a,W_u],L_{V_a}\xi^u-L_{W_u}\nu^a+d(\nu^a(W_u))) $$
$$=([V_a,W_u],L_{V_a}\xi^u - i(W_u)d\nu^a) = (L_{V_a}W_u,L_{V_a}\xi^u)$$
$$- (0,(W_u\alpha^a_b)dx^b)=f^a\mathcal{X}_a+\varphi^u\Xi_u,$$
where we have used the local decomposability condition $\alpha^{'a}_h=0$.
Since $[V_a,W_u]$ does not contain $\partial/\partial x^a$, we must
have $f^a=0$ and  (\ref{auxZ}) follows.\end{proof}

Now, we can prove
\begin{prop}\label{strtransv} Let $E$ be a strongly, locally decomposable,
integrable, big-isotropic structure on a manifold $M$. Then, the
transversal structure $E^{tr}_{x_0}$
is well defined up to a structure preserving diffeomorphism.\end{prop}
\begin{proof} The proof given for the Dirac structures in \cite{DW} also holds here.
First we look at local transversal submanifolds $\mathcal{T}_0,\mathcal{T}_1$
of the characteristic slice $\mathcal{U}$ at $x_0\neq x_1\in\mathcal{U}$.
Then, there exist
a tubular neighborhood such that $\mathcal{T}_0,\mathcal{T}_1$ are
fibers of the tubular foliation $\mathcal{F}$ and a diffeomorphism $\Phi$,
which is a composition of transformations of flows of $\mathcal{F}$-projectable
vector fields $Z\in H(E,\mathcal{F})$, that sends $\mathcal{T}_0$
onto $\mathcal{T}_1$. By Lemma \ref{lemaptisotopie}, $\Phi$ sends
the transversal structure on $\mathcal{T}_0$ onto the transversal
structure on $\mathcal{T}_1$. Now, if we have two different local
transversal submanifolds $\mathcal{T}_0,\mathcal{T}'_0$ at the
same point $x_0$ of $\mathcal{U}$, we take a loop of $\mathcal{U}$
at $x_0$, break it into a finite number of pieces, and go from
$\mathcal{T}_0$ to $\mathcal{T}'_0$ through intermediate
transversal submanifolds defined at the breaking points.
The composition of the diffemorphisms $\Phi$ defined as above between
the intermediate transversal manifolds gives us the required
equivalence of the transversal structures on $\mathcal{T}_0$ and
$\mathcal{T}'_0$.\end{proof}
\section{Reduction of big-isotropic structures}
Recently, the generalization of symplectic reduction to Dirac
structures was discussed by several authors, in particular
\cite{{BG},{SXu}}. In this section we discuss a reduction scheme
for big-isotropic structures.

We begin by defining a push forward procedure. The notation will
be similar to that of Proposition \ref{proppull}; in particular,
we consider the mapping $f:N\rightarrow M$ and the corresponding
points $y=f(x)$. But, we start with a rank $k$ subbundle
$E\subseteq T^{big}N$. Then, we define the push forward of $E$ by
\begin{equation}\label{defpush} f_*(E_x)=
\{(f_*X,\alpha)\,/\,X\in T_xN,\,\alpha\in T^*_yM,\,(X,f^*_y\alpha)\in E_x\}.
\end{equation}
\begin{prop}\label{proppush} If the bundle $E$ is a big-isotropic
structure then $f_*(E_x)$ is an isotropic subspace of
$(T^{big}_yM,g_M)$ and its orthogonal space is $f_*(E'_x)$.
Furthermore, iff
\begin{equation}\label{condpush} E_x\cap ker\,f_{*x}= E'_x\cap
ker\,f_{*x}\end{equation} the dimension of $f_*(E_x)$ is $m-n+k$.
\end{prop} \begin{proof} The isotropy of $f_*(E_x)$ is obvious.
Then, let us define the space
\begin{equation}\label{eqSigma} \Sigma=\{(X,f^*\alpha)\in E_x\,/\,
X\in T_xN,\,\alpha\in T_y^*M\}=E_x\cap(T_xN\oplus
im\,f^*_y)\end{equation} $$=(E'_x+ker\,f_{*x})^{\perp_{g_N}}.$$
The last equality (\ref{eqSigma}) follows from
\begin{equation}\label{xplicSigma} (ker\,f_{*x})^{\perp_{g_N}} =
T_xN\oplus im\,f^*_y,\end{equation} which holds because the right
hand side is included
in the left hand side and the former has the dimension required
for the orthogonal space of the latter. The correspondence
$(X,f^*\alpha)\mapsto (f_*X,\alpha)$ produces an isomorphism
\begin{equation}\label{isomSigma} \Sigma/(E_x\cap ker \,f_{*x})
\approx f_*(E_x)/ker\,f^*_y, \end{equation} whence
$$dim\,f_*(E_x)=dim\,ker\,f^*_y+dim\,\Sigma-dim(E_x\cap
ker\,f_{*x})$$ $$\stackrel{(\ref{eqauxptf})}{=}
m-rank_xf+dim\,\Sigma-2n+dim(E_x\cap ker\,f_{*x})^{\perp_{g_N}}$$
$$\stackrel{(\ref{xplicSigma})}{=} m-rank_xf+dim\,\Sigma-2n+dim
((T_xN\oplus im\,f^*_y)+E'_x)$$ $$=m-rank_xf+dim\,\Sigma-2n+dim
(T_xN\oplus im\,f^*_y)$$ $$+dim\,E'_x - dim((T_xN\oplus im\,f^*_y)
\cap E'_x)$$ $$=m+n-k+(dim\,\Sigma-dim\,\Sigma'),$$ where $\Sigma'$
is the space (\ref{eqSigma}) for $E'$.

The same calculations with the roles of $E$ and $E'$ interchanged
give $$dim\,f_*(E'_x)=m-n+k-(dim\,\Sigma-dim\,\Sigma'),$$ hence,
$dim\,f_*(E_x)+dim\,f_*(E'_x)=2m$. Since it is trivial to check
that the two spaces are $g_M$-orthogonal, we get the required
relation $f_*(E'_x)\perp_{g_M}f_*(E_x)$.

Furthermore, from the last expression of $\Sigma$ in
(\ref{eqSigma}) and the similar expression of $\Sigma'$ we get
$$dim\,\Sigma -dim\,\Sigma'=dim\,\Sigma^{'\perp_{g_N}}
-dim\,\Sigma^{\perp_{g_N}}$$ $$= dim\,E_x+dim\,ker\,f_{*x}
-dim(E_x\cap ker\,f_{*x})$$
$$-dim\,E'_x-dim\,ker\,f_{*x}+dim(E'_x\cap ker\,f_{*x})$$
$$=2k-2n +(dim(E'_x\cap ker\,f_{*x}) - dim(E_x\cap ker\,f_{*x})).$$
Accordingly, we obtain \begin{equation}\label{finaldimpush}
dim\,f_*(E_x)=m-n+k +(dim(E'_x\cap ker\,f_{*x}) - dim(E_x\cap
ker\,f_{*x})),\end{equation} which justifies the last assertion of
the proposition.\end{proof}
\begin{rem}\label{obscodim2} {\rm  Proposition \ref{proppush} holds for any
linear mapping $l:V\rightarrow W$ between linear spaces and any
isotropic subspace $E\subseteq V\oplus V^*$. The
significance of the equality $dim\,f^*(E_y)=m-n+k$ is the
codimensional invariance property $n-dim\,E_x=m-dim\,f_*(E_x)$.
Hypothesis
(\ref{condpush}) holds in the Dirac case $(E=E')$ and in the case
where $f$ is an immersion $(kef\,f_*=0)$. If $f$ is a submersion
and the difference between the dimensions of $\Sigma,\Sigma'$ is
constant then $dim\,f_*(E_x)$ is the
same $\forall x\in N$. However, there may not be a well defined
subbundle $f_*E$ because $f^-1(y)$ may have more than one
point.}\end{rem}
\begin{prop}\label{operatiiinverse} If $f$ is a submersion then
$f_*f^*(E_y)= E_y$. If $f$ is an immersion then $f^*f_*(E_x)= E_x$.
\end{prop} \begin{proof} By looking at the formulas (\ref{pullback}),
(\ref{defpush}) we see that $(Z,\lambda)\in f_*f^*(E_y)$ iff
$(Z,\lambda)= (f_*X,\alpha+\beta)$ where $(f_*X,\alpha)\in E_y$ and
$\beta\in ker\,f^*_y$. Since, if $f$ is a submersion $ker\,
f^*_y=0$, we get $(Z,\lambda)\in E_y$, hence, $f_*f^*(E_y)\subseteq
E_y$. On the other hand, again since $f$ is a submersion, any pair
of $E_y$ is of the form $(f_*X,\alpha)$, hence, of the form
$(Z,\lambda)\in f_*f^*(E_y)$. Thus, $E_y\subseteq f_*f^*(E_y)$ and
the first assertion is proven. The proof of the second assertion is
similar.
\end{proof}

Another kind of preparation that we need concerns the notion of
projectability. Let $M$ be an $m$-dimensional, differentiable
manifold and $\mathcal{F}$ a foliation of $M$ by $p$-dimensional
leaves. In what follows, the terms projectable and foliated are
synonymous and describe objects related with corresponding objects
of the space of leaves via the natural projection. In particular, a
vector field is foliated if it can be projected onto the space of
leaves and a differential form is projectable if it is the pullback
of a form on the space of leaves (such forms are called basic forms
by most authors).
\begin{defin}\label{defproj} {\rm An arbitrary subbundle
$E\subseteq T^{big}M$ of rank $k$ is called {\it foliated} or {\it
projectable} if each point $x\in M$ has an open neighborhood $U$
such that the quotient manifold $Q_U=U/ (\mathcal{F}|_U)$ is endowed
with a subbundle $\Delta_U\subseteq T^{big}Q_U$ that satisfies the
condition $E|_U=\pi^*(\Delta_U)$ $(\pi:U\rightarrow U/
(\mathcal{F}|_U))$.}\end{defin}

In Definition \ref{defproj}, $\pi$ is the natural projection and
$\pi^*(\Delta_U)$ is obtained by the pullback of $\Delta_U$, i.e.,
\begin{equation}\label{backD} E_x=\{(X,\pi^*\alpha)\,/\,X\in T_xM,\,
(\pi_*X,\alpha)\in\Delta_{\pi(x)}\}\hspace{2mm}(x\in
U).\end{equation}
\begin{prop}\label{conddeproj} The subbundle $E$ is foliated iff
it satisfies the following two conditions: a) $E\supseteq
T\mathcal{F}$, b) each point $x\in M$ has an open neighborhood $U$
such that $\Gamma E|_U$ has a basis $(X_i,\xi^i)$ $(i=1,...,k)$ with
projectable vector fields $X_i$ and projectable $1$-forms $\xi^i$.
Furthermore, condition b) may be replaced by b') every $Y\in\Gamma
T\mathcal{F}$ is an infinitesimal automorphism of $E$, i.e.,
$\forall (X,\xi)\in\Gamma E$ one has $(L_YX, L_Y\xi)\in\Gamma
E$.\end{prop}
\begin{proof} If $E$ is projectable, since $(0,0)\in\Delta_U$,
formula (\ref{backD}) shows that a) holds. Then, assume that
$(V_u,\lambda^u)$ $(u=1,...,q=m-p)$ is a local basis of the cross
sections of $\Delta_U$ and put
$V_u=[X_u]_{T\mathcal{F}},\xi_u=\pi^*\lambda^u$ $(X\in\chi^1(U))$,
where the brackets denote equivalence classes modulo $T\mathcal{F}$.
It follows easily that $\Gamma E|_U$ has a local basis that consists
of $(X_u,\xi^u)$ and of a basis of $\Gamma T\mathcal{F}$. Hence b)
also holds. Conversely, if we have the properties a), b), we can
change the basis provided by b) to a basis of the form
$((X_u,\xi^u),(Y_a,0))$ where $Y_a$ is a local basis of
$T\mathcal{F}$ and $([X_u]_{T\mathcal{F}},\lambda^u)$,
$\xi^u=\pi^*\lambda^u$, is a basis for the local structure
$\Delta_U$ required by Definition \ref{defproj}. For the last
assertion, if b) holds and if we put
$(X,\xi)=\sum_{i=1}^kf^i(X_i,\xi^i)$ then
$$(L_YX,L_Y\xi)=\sum_{i=1}^k(Yf^i)(X_i,\xi^i)\in\Gamma E.$$
Conversely, assume that a), b') hold. For any $x\in M$, we can
take a cubical, open neighborhood with $ \mathcal{F}$-adapted
local coordinates $(z^a,y^u)$ (i.e., $\mathcal{F}$ has the local
equations $dy^u=0$ and $z^a(x)=0,y^u(x)=0$). Then, we can take an
arbitrary basis of $E$ of the form $(\partial/\partial
z^a,0),(V_u,\lambda^u))$ along the transversal slice $z^a=0$ and
move $(V_u,\lambda^u)$ along the slice $y^u=0$ by the linear
holonomy of the foliation $ \mathcal{F}$ (e.g., \cite{Mol}). By
hypothesis b') the flows of the tangent vectors of the leaves
preserve $E$, hence, the result of the previous procedure consists
of cross sections of $E$ and we get a projectable basis as
required by b).\end{proof}

We may apply the previous general results to speak of projectable,
big-isotropic structures, integrable or not, and of projectable
(almost) Dirac structures. First, we notice the following
corollary of Proposition \ref{conddeproj}.
\begin{corol}\label{integrdaproj} An integrable,
big-isotropic structure $E$ of $M$ is foliated iff $E\supseteq
T\mathcal{F}$.
\end{corol} \begin{proof}. For a foliated subbundle $E$ we have
$E\supseteq T\mathcal{F}$ by a) of Proposition \ref{conddeproj}.
Conversely, since $\Gamma E$ is closed by Courant brackets, with
the notation of b'), Proposition
\ref{conddeproj}, we have $$[(Y,0),(X,\xi)]=
(L_YX,L_Y\xi)\in \Gamma E.$$ Hence, conditions a) and b') hold and
we are done.\end{proof}
\begin{prop}\label{integrproj} For a big-isotropic, projectable subbundle $E\subseteq
T^{big}M$, $\Gamma E$ is closed by Courant brackets iff the spaces
$\Gamma\Delta_U$ of the local projected subbundles $\Delta_U$ are
closed by Courant brackets. \end{prop}
\begin{proof} For two differentiable, local, cross sections of $E$ of the
form prescribed by (\ref{backD}) the Courant bracket has the following expression
\begin{equation}\label{crCr1}
[(X,\pi^*\alpha),(Y,\pi^*\beta)]=([X,Y],\pi^*(L_{\pi_*X}\beta -
L_{\pi_*Y}\alpha)\end{equation}
$$+\frac{1}{2}d(\alpha(\pi_*Y)-\beta(\pi_*X))) $$
and it is related by (\ref{backD}) with the Courant bracket $[(\pi_*X,\alpha),
(\pi_*Y,\beta)]$ on $Q_U$. Thus, if $\Gamma\Delta$ is closed by Courant
brackets, the left hand side of (\ref{crCr1}) is in $\Gamma E$. Furthermore,
the Courant bracket
satisfies the property that, $\forall f\in C^\infty(M)$, one has
\begin{equation}\label{propdealg} [(X,\alpha),f(Y,\beta)] =
f[(X,\alpha),(Y,\beta)] +(Xf)(Y,\beta)
\end{equation} $$ - g((X,\alpha),f(Y,\beta))(0,df)$$
and the isotropy of $E$ yields
$$[f(X,\pi^*\alpha),h(Y,\pi^*\beta)]=
fh[(X,\pi^*\alpha),(Y,\pi^*\beta)]$$
$$+f(Xh)(X,\pi^*\alpha) -h(Yf)(Y,\pi^*\beta),$$
where $f,h\in C^\infty(M)$ may not be projectable. Accordingly, the
existence of projectable bases of $E$ shows that the bracket-closure
of $\Gamma\Delta_U$ implies the bracket-closure of $\Gamma E$. The
converse result is a straightforward consequence of the existence of
the projectable bases of $E$.
\end{proof} \begin{rem}\label{obsTstea} {\rm From (\ref{propdealg}) it follows
that  if $E\subseteq T^{big}M$ is closed by Courant brackets then either
$E$ is isotropic
or $E\supseteq T^*M$.}\end{rem}

\begin{prop}\label{projE'} If $E$ is a projectable, big-isotropic
structure on the foliated manifold $(M,\mathcal{F})$ its orthogonal
bundle $E'$ is projectable as well. Moreover, the corresponding,
local, projected bundles $\Delta_U$ are big-isotropic on the local
transversal manifolds $\mathcal{Q}_U$, the orthogonal bundles
$\Delta'_U$ are the local, projected bundles of $E'$, and $E$ is
integrable iff the structures $\Delta_U$ are integrable.\end{prop}
\begin{proof} It is easy to check conditions a) and b') of
Proposition \ref{conddeproj} for $E'$. (In particular, for b')
take $(X,\alpha)\in\Gamma E,(Y,\beta)\in\Gamma E', Z\in
T\mathcal{F}$, express $Z(g((X,\alpha),(Y,\beta)))=0$ and use b')
for $E$.) For the other assertions, use Proposition
\ref{integrproj}.\end{proof}
\begin{example}\label{exstrhamilt} {\rm Let $P\in\chi^2(M)$
be a bivector field on $(M,\mathcal{F})$ and define
\begin{equation}\label{eqDhamilt}
D_P=T\mathcal{F}\oplus\{(\sharp_P\alpha,\alpha)\,/\,\alpha\in
ann(T\mathcal{F})\}.\end{equation} This is an almost Dirac
structure, which does not depend on the $T\mathcal{F}$-component
of $P$ in the following sense. If $\nu\mathcal{F}$ is a normal
bundle of $
\mathcal{F}$, i.e., $TM=\nu\mathcal{F}\oplus T\mathcal{F}$, and if
we put
\begin{equation}\label{Pcucomp} P= \frac{1}{2}P^{ab}Z_a\wedge Z_b
+P^{au}Z_a\wedge Y_u +\frac{1}{2}P^{uv}Y_u\wedge
Y_v,\end{equation} where $(z^a,y^u)$ are the local coordinates
used in the proof of Proposition \ref{conddeproj} and
$Z_a=\partial/\partial z^a,Y_u=\partial/\partial y^u-t^a_uZ_a$ (for some
coefficients $t^a_u$),
$D_P$ is spanned by $(Z_a,0)$ and $(P^{uv}Y_v,dy^u)$. Conditions a), b')
of Proposition \ref{conddeproj} show that the projectability of $D_P$ is
equivalent with the
projectability of the bivector field $P$, i.e.,
$\partial P^{uv}/\partial z^a=0$. Furthermore, using Proposition \ref{projE'},
we see that the integrability
of $D_P$ holds iff $[P,P]|_{ann(T\mathcal{F})}=0$ (Schouten-Nijenhuis bracket),
hence $D_P$
is equivalent with the {\it transversely Hamiltonian structure}
defined by $P$ on $(M,\mathcal{F})$ \cite{IV-Lux}.}
\end{example}
\begin{example}\label{exFpres} {\rm In a similar
way, let $\omega$ be a foliated $2$-form on $(M,\mathcal{F})$ and
$\nu\mathcal{F}$ be a chosen normal bundle. Define
\begin{equation}\label{defFpres} D_\omega= T\mathcal{F}\oplus
\{(X,\flat_\omega X)\,/\,X\in \nu\mathcal{F}\}.\end{equation} Then
$D_\omega$ is an almost Dirac structure. If we put
\begin{equation}\label{compomega} \omega=
\frac{1}{2}\omega_{uv}(y)dy^u\wedge dy^v,\end{equation} we see
that $D_\omega$ is spanned by
$(Z_a,0),(Y_u,\omega_{uv}Y_v)$, therefore, $D_\omega$ is
foliated and independent on the choice of $\nu\mathcal{F}$.
Furthermore, $D_\omega$ is integrable iff $\omega$ is
closed, i.e., $D_\omega$ is equivalent with a foliated presymplectic form.}
\end{example}

Now, we will discuss the concept of reduction. The general geometric
framework may be described as follows. Let $M$ be a manifold, $E$ a
subbundle of $T^{big}M$ and $\iota:N\hookrightarrow M$ an embedded
submanifold. Assume that $N$ is $E$-{\it proper}, meaning that the
pullback $\iota^*E$ is differentiable. Then, assume that $N$ is
endowed with a foliation $\mathcal{F}$ such that
$T\mathcal{F}\subseteq\iota^*E$, a condition that is equivalent with\vspace{2mm}\\
(r)\hspace{2cm} $\forall Z\in T\mathcal{F}$, $\exists\alpha\in ann\,TN$ such that
$(Z,\alpha)\in E|_N$\vspace{2mm}\\
and will be called the {\it reducibility condition}. Then, $\iota^*E$ may be foliated.
If it is so and if the quotient space $N/\mathcal{F}$ is a
paracompact, Hausdorff manifold $Q$, $\iota^*E$ projects to a
subbundle $E_N^{red}\subseteq T^{big}Q$ that has the restrictions
$\Delta_U$ of Definition  \ref{defproj} over the projections $\pi(U)$.
Accordingly, $\pi^* E^{red}_N=\iota^*E$ and Proposition \ref{operatiiinverse}
shows that we may write
\begin{equation}\label{redX} E_N^{red}=\pi_*(\iota^*E).
\end{equation}
The bundle $E_N^{red}$ will be called the {\it reduction of
$E$ via $(N,\mathcal{F})$}. This is a generalization of the
framework described in the paper of Sti\'enon-Xu \cite{SXu}.

Now, the following result is immediate.
\begin{prop}\label{propdered} Let $E$ be an integrable,
big-isotropic structure on the manifold $M$. Let $\iota:N\rightarrow
M$ be an embedded submanifold such that the fields of subspaces
$E\cap(TN\oplus T^*_NM),E'\cap(TN\oplus T^*_NM)$ are of a constant
dimension. Let $\mathcal{F}$ be a foliation of $N$ by the fibers of
a submersion $\pi:N\rightarrow Q$ such that the natural projection
\begin{equation}\label{2red}
pr_{T\mathcal{F}}:E\cap(T\mathcal{F}\oplus ann\,TN)
\rightarrow T\mathcal{F}\end{equation} is
surjective. Then, there exists a well defined, integrable, reduced,
big-isotropic structure $E^{red}$ on $Q$ that satisfies the
condition {\rm(\ref{redX})}.\end{prop} \begin{proof} By Corollary
\ref{cazdedif}  $\iota^*E$ is an integrable, big-isotropic
structure of $N$ and condition (\ref{2red}) is equivalent with the reducibility condition (r).  Thus, the reduced structure $E^{red}$ may exist and
Corollary \ref{integrdaproj} tells us that $E^{red}$ exists indeed. Then, Proposition
\ref{projE'} shows that $E^{red}$ is an integrable, big-isotropic structure on $Q$.\end{proof}
\begin{corol}\label{Gred} Let $E$ be an integrable,
big-isotropic structure on the manifold $M$. Assume that the connected, Lie
group $G$ acts on $M$ and the action preserves $E$ and fixes an embedded
submanifold $\iota:N\rightarrow M$. Assume that the restriction of
the action of $G$ to $N$ is proper and free and denote by
$\mathcal{F}$ the foliation of $N$ by the orbits of $G$. Finally,
assume that the reducibility condition {\rm(r)} holds for the
infinitesimal transformations $Z$ of $G$ on $N$. Then, there
exists a Hausdorff manifold $Q=N/G$ endowed with a reduced,
integrable, big-isotropic structure $E^{red}$.\end{corol}
\begin{proof} Under the hypotheses, condition (r) holds for any
$Z\in\chi^1(N)$ and the fields of subspaces
$E\cap(TN\oplus T^*_NM)$, $E'\cap(TN\oplus T^*_NM)$ have a
constant dimension.\end{proof}
\begin{example}\label{Pred} {\rm We apply Proposition
\ref{propdered} to a Poisson structure $E_P=E'_P=graph\,\sharp_P$,
where $P\in\chi^2(M)$ is a Poisson bivector field on $M$. Assume
that the submanifold $\iota^*:N\hookrightarrow M$ is such that
$dim(E_P\cap(TN\oplus T^*_NM))=const.$ and that $N$ has a foliation
$\mathcal{F}$ with the quotient manifold $Q=M/\mathcal{F}$. Then,
$N$ has the Dirac structure
$$\iota^*E_P=\{(\sharp_P\alpha,\iota^*\alpha)\,/\,\sharp_P\alpha\in TN,
\alpha\in T^*_NM\}.$$ Furthermore, the reducibility condition (r) is
equivalent to \begin{equation}\label{2redinex}
T\mathcal{F}\subseteq\sharp_P(ann\, TN).\end{equation} Therefore, if
we assume that (\ref{2redinex}) holds, $Q$ has the reduced Dirac structure
$E^{red}=\pi_*\iota^*(E_P)$, which, pointwisely, turns out to be
\begin{equation}\label{redinex}
E^{red}=\{(\pi_*\sharp_P\tilde{\lambda},\lambda)\,/\, \lambda\in
T^*Q,\,[\tilde{\lambda}]_{ann\,TN}=\pi^*\lambda,\,
\sharp_P\tilde{\lambda}\in TN\}.\end{equation} The reduced structure
$E^{red}$ is Dirac. If we want a Poisson reduced structure we have
to add the condition $E^{red}\cap TQ=0$, which, by (\ref{redinex}),
is equivalent to
$$\sharp_P(ann\,TN)\cap TN\subseteq T\mathcal{F}.$$
A more general Poisson reduction scheme, where $T\mathcal{F}=V\cap TN$
for a vector bundle $V$ over $N$ that satisfies adequate conditions, was given by
Marsden and Ratiu \cite{MR}. The present example corresponds to the
case $V=\sharp_P(ann\,TN)$, while, however, we do not ask $\sharp_P(ann\,TN)$
to be a regular vector bundle.}
\end{example} \vspace*{2mm}
{\it Acknowledgement}. Part of the work on this paper was done
during the author's visit to the Bernoulli Center of the \'Ecole
Polytechnique F\'ed\'erale de Lausanne, Switzerland, in
June-August 2006, and the author wishes to express his gratitude
to the Center and to professor Tudor Ratiu, the director of the
Bernoulli Center, in particular, for the invitation and support.
\hspace*{7.5cm}{\small \begin{tabular}{l} Department of
Mathematics\\ University of Haifa, Israel\\ E-mail:
vaisman@math.haifa.ac.il \end{tabular}}
\end{document}